\documentclass[10pt,twoside,leqno]{article}
\pdfoutput=1
\usepackage[utf8]{inputenc}
\usepackage[T1]{fontenc}
\usepackage[full]{textcomp}

\usepackage{csquotes}
\usepackage[english]{babel}

\usepackage{mathtools}
\mathtoolsset{mathic} 

\usepackage{cfr-lm}
\usepackage{dsfont}  

\usepackage{microtype}
\usepackage[
nohead,nomarginpar,
paper=a4paper,
]{geometry}

\usepackage{lastpage}
\usepackage{fancyhdr}
\usepackage{longtable}
\usepackage{booktabs}
\usepackage{tabu}

\usepackage[inline]{enumitem}
\setlist{noitemsep,nosep,listparindent=\parindent}
\setlist[itemize]{label=\guillemotright}
\setlist[enumerate,1]{ref=\thesubsection.\arabic*}
\setlist[enumerate,2]{label=\alph*.,ref=\theenumi.\alph*}

\usepackage{cjhebrew}

\usepackage[backend=biber,doi=false,url=false,isbn=false,%
safeinputenc,style=quasialphabetic,citestyle=alphabetic]{biblatex}
\bibliography{\jobname.bib}
\defbibheading{bibliography}[\bibname]{\sectionstar{#1}}

\usepackage{ifthen}
\makeatletter
\renewcommand{\part}[1]{%
 \cleardoublepage%
 \vbox{\null\vskip90pt%
 \normalfont\fontsize{20pt}{30pt}\selectfont%
 \baselineskip=30pt%
 \scshape\noindent\textls*{#1}\par}%
 \addcontentsline{toc}{part}{#1}%
 \@afterindentfalse%
 \@afterheading%
}
\renewcommand{\section}[1]{%
 \vskip2\baselineskip\penalty-250%
 \refstepcounter{section}%
 \vbox{\normalfont\fontsize{12pt}{15pt}\selectfont%
  \centering\scshape\noindent\textls*{\thesection\quad#1}%
  \par}
 \nobreak
 \addcontentsline{toc}{section}{\protect\numberline{\thesection} #1}%
 \@afterindentfalse%
 \@afterheading%
} 
\newcommand{\sectionstar}[1]{%
 \vskip2\baselineskip\penalty-250
 \vbox{\normalfont\fontsize{12pt}{15pt}\selectfont%
  \centering\scshape\noindent\textls*{#1}%
  \par}
 \nobreak\vskip15pt
 \@afterindentfalse%
 \@afterheading%
} 
\renewcommand{\paragraph}[1]{\par\bigskip\refstepcounter{subsection}%
 \global\setbox\@labels\hbox{%
  \mbox{\normalfont\normalsize\scshape\noindent\thesubsection%
   \ifthenelse{\equal{#1}{}}%
   {}%
   {\ \textls{#1.}}%
   \ ---}%
 }%
 \global\@inlabeltrue
 \everypar{%
  \if@inlabel
  \global\@inlabelfalse
  {\setbox\z@ \lastbox}%
  \unhbox\@labels
  \spacefactor \sfcode `.\relax
  \space
  \penalty\z@
  \fi
 }
 \ignorespaces
}
\newcommand{\readmetitle}{\par\vskip\baselineskip%
 {\normalfont\normalsize\scshape\noindent%
  \textls{Readme.}\ ---}
}
\newenvironment{readme}{\readmetitle\itshape}{\normalfont}
\renewcommand\tableofcontents{%
 \sectionstar{\contentsname}%
 \@starttoc{toc}%
}
\renewcommand*\l@part[2]{%
 \addvspace{15pt \@plus\p@}%
 \noindent{\leavevmode%
  \scshape\textls{#1\qquad#2}%
 }\par\nobreak%
}
\renewcommand*\l@section[2]{%
 \setlength\@tempdima{\parindent}%
 \noindent
 {\leavevmode%
  \hskip\parindent#1\qquad#2%
 }\par\nobreak%
}
\makeatother

\usepackage{amsmath,amssymb}  
\usepackage{mathrsfs}         
\usepackage{mathabx}

\usepackage[thmmarks,amsmath]{ntheorem}
\usepackage{hyperref}
\usepackage{thmtools}

\numberwithin{equation}{subsection}

\declaretheoremstyle[headformat=swapnumber,headpunct={.\ ---},%
headfont=\normalfont\scshape\lsstyle,bodyfont=\itshape,%
spaceabove=0pt,spacebelow=0pt,%
preheadhook={\bigskip}]{theorem}
\declaretheorem[style=theorem,sibling=subsection]{theorem}
\declaretheorem[style=theorem,sibling=subsection]{proposition}
\declaretheorem[style=theorem,sibling=subsection]{lemma}

\declaretheorem[style=theorem,sibling=subsection]{conjecture}

\declaretheoremstyle[headformat=swapnumber,headpunct={.\ ---},%
headfont=\normalfont\scshape\lsstyle,bodyfont=\normalfont,%
spaceabove=0pt,spacebelow=0pt,%
preheadhook={\bigskip}]{definition}
\declaretheorem[style=definition,sibling=subsection]{definition}

\declaretheorem[style=definition,sibling=subsection]{example}
\declaretheorem[style=definition,sibling=subsection]{remark}
\declaretheorem[style=definition,sibling=subsection]{notation}

\def\qedsquare{\ \boxvoid}

\declaretheoremstyle[headpunct={\!.},headfont=\itshape,bodyfont=\normalfont,%
qed=\ensuremath{\qedsquare},spaceabove=0pt,spacebelow=0pt]{proof}
\declaretheoremstyle[headpunct={\!.},headfont=\itshape,bodyfont=\normalfont,%
qed=\ensuremath{\qedsquare},spaceabove=0pt,spacebelow=0pt]{nonumberproof}
\declaretheorem[style=proof,numbered=no]{proof}

\declaretheoremstyle[headformat=swapnumber,headpunct={.\ ---},%
headfont=\itshape,bodyfont=\normalfont,qed=\ensuremath{\qedsquare},%
spaceabove=0pt,spacebelow=0pt,%
preheadhook={\bigskip}]{nproof}

\usepackage{cleveref}
\crefname{condition}{condition}{conditions}
\crefname{conjecture}{conjecture}{conjectures}
\crefname{construction}{construction}{constructions}
\crefname{corollary}{corollary}{corollaries}
\crefname{diagram}{diagram}{diagrams}
\crefformat{subsection}{\S#2#1#3}
\crefformat{enumi}{\S#2#1#3}
\crefformat{nproof}{\S#2#1#3}
\creflabelformat{equation}{#2#1#3}

\newcommand{\id}{\textnormal{id}}
\newcommand{\ev}{\textnormal{ev}}

\newcommand{\into}{\hookrightarrow}
\newcommand{\onto}{\twoheadrightarrow}
\newcommand{\longto}{\longrightarrow}

\newcommand{\Hom}{\textnormal{Hom}}
\newcommand{\End}{\textnormal{End}}
\newcommand{\Isom}{\textnormal{Isom}}

\newcommand{\Aut}{\textnormal{Aut}}

\newcommand{\iAut}{\underline{\Aut}}

\newcommand{\Mat}{\textnormal{Mat}}

\newcommand{\NN}{\mathbb{N}}
\newcommand{\ZZ}{\mathbb{Z}}
\newcommand{\QQ}{\mathbb{Q}}
\newcommand{\QQbar}{\bar{\QQ}}
\newcommand{\QQl}{\QQ_{\ell}}

\newcommand{\QQp}{\QQ_{p}}

\newcommand{\BB}{\mathrm{B}}
\newcommand{\RR}{\mathbb{R}}
\newcommand{\CC}{\mathbb{C}}
\newcommand{\HQ}{\mathbb{H}}

\newcommand{\primes}{\mathscr{L}}

\newcommand{\Spec}{\textnormal{Spec}}

\newcommand{\DelS}{\mathbb{S}}

\newcommand{\Gal}{\textnormal{Gal}}

\newcommand{\HdR}{\HH_{\dR}}

\newcommand{\HT}{\textnormal{HT}}
\newcommand{\BHT}[1]{\textnormal{B}_{\HT,#1}}
\newcommand{\gr}{\textnormal{gr}}
\newcommand{\GdR}{\GG_{\dR}}

\newcommand{\Vect}{\textnormal{Vect}}
\newcommand{\grVect}{\textnormal{grVect}}

\newcommand{\Rep}{\textnormal{Rep}}

\newcommand{\FpHS}{\textnormal{FpHS}}

\makeatletter
\def\cpwith[#1]#2{\textnormal{c.p.}_{#1}(#2)}
\def\cpwithout#1{\textnormal{c.p.}(#1)}
\def\cp{\@ifnextchar[{\cpwith}{\cpwithout}}
\makeatother
\makeatletter
\def\Gmwith[#1]{\mathbb{G}_{\textnormal{m},#1}}
\def\Gmwithout{\mathbb{G}_{\textnormal{m}}}
\def\Gm{\@ifnextchar[{\Gmwith}{\Gmwithout}}
\makeatother
\newcommand{\GL}{\textnormal{GL}}

\newcommand{\Lie}{\textnormal{Lie}}

\newcommand{\ab}{\textnormal{ab}}
\newcommand{\der}{\textnormal{der}}
\newcommand{\ad}{\textnormal{ad}}
\newcommand{\ha}{\textnormal{ha}}

\newcommand{\Fib}{\textnormal{Fib}}

\newcommand{\Proj}{\textnormal{Proj}}

\newcommand{\cl}{\textnormal{cl}}

\newcommand{\dR}{\textnormal{dR}}
\newcommand{\et}{\textnormal{\'{e}t}}
\newcommand{\sing}{\textnormal{sing}}

\newcommand{\HH}{\textnormal{H}}
\newcommand{\Hl}{\HH_{\ell}}
\newcommand{\Hp}{\HH_{p}}
\newcommand{\Hlambda}{\HH_{\lambda}}
\newcommand{\HB}{\HH_{\textnormal{B}}}
\newcommand{\HLambda}{\HH_{\Lambda}}

\newcommand{\Mot}{\textnormal{Mot}}

\newcommand{\Zentrum}{\textnormal{Z}}
\newcommand{\GG}{\textnormal{G}}
\newcommand{\GB}{\GG_{\textnormal{B}}}
\newcommand{\Gp}{\GG_{p}}
\newcommand{\Gpc}{\Gp^{\circ}}
\newcommand{\Gl}{\GG_{\ell}}
\newcommand{\Glc}{\Gl^{\circ}}
\newcommand{\Glambda}{\GG_{\lambda}}
\newcommand{\Glambdac}{\Glambda^{\circ}}

\newcommand{\HS}{\textnormal{HS}}

\newcommand{\Res}{\textnormal{Res}}

\newcommand{\chrc}{\textnormal{char}}

\newcommand{\Tangen}[1]{\langle #1 \rangle^{\otimes}}

\newcommand{\cm}{\textsc{cm}}

\newcommand{\MTC}{\textnormal{MTC}}

\newcommand{\DtA}{\textnormal{A}}
\newcommand{\DtB}{\textnormal{B}}
\newcommand{\DtC}{\textnormal{C}}
\newcommand{\DtD}{\textnormal{D}}
\newcommand{\DtE}{\textnormal{E}}
\newcommand{\DtF}{\textnormal{F}}
\newcommand{\DtG}{\textnormal{G}}

\usepackage{tikz}
\usetikzlibrary{calc,
	arrows.meta,
	cd,
	decorations.markings,
	decorations.pathmorphing,
	positioning,
	positioning,
	shadows,
	shapes,
	shapes.geometric
}
\tikzset{Dynkin/.style = {
		every node/.style = {circle, draw, semithick, inner sep=2pt, fill=white,
			copy shadow={fill=black, shadow xshift=0.3pt, shadow yshift=-0.1pt},
		},
		every label/.append style = {label distance=-2.5pt,
   rectangle,draw=none,font=\footnotesize,
   inner sep=1ex,text depth=1pt,
			every shadow/.style={opacity=0}
		},
  special/.style={rectangle, inner sep=2.4pt},
  symplectic/.style={correct forbidden sign},
		a0/.style = {fill=gray!50},
  diagramlabel/.style = {draw=none,opacity=0},
		doubledynkin/.style={double distance=2pt,
			decoration={markings,
				mark=at position 0.6 with {\arrow[semithick]{Straight Barb[length=5pt]}},
			},
			postaction={decorate},
		}
	}
}

\def\adef#1{\catcode`#1=13 \bgroup \lccode`\~=`#1\lowercase{\egroup\def~}}
\long\def\addto#1#2{\expandafter\def\expandafter#1\expandafter{#1#2}}
\def\activebraces#1#2{\adef#1{\ifmmode#1\else{\textup#1}\fi}\adef#2{\ifmmode#2\else{\/\textup#2}\fi}}
\let\emphOri=\emph
\def\emph{\bgroup\activebraces()\activebraces[]\emphA}
\def\emphA#1{\emphOri{#1}\egroup}
\addto\itshape{\activebraces()}

\def\title{The Mumford--Tate conjecture for products of abelian varieties}
\def\author{Johan Commelin}

\usepackage{datetime}
\def\date{\dayofweekname{\day}{\month}{\year},
 the \ordinaldate{\day} of \monthname, \number\year}
\begin{document}
\begin{center}\Large\scshape
\textls*{\title}
\end{center}

\medskip

\noindent\textit{by} \quad \author \hfill \date

\bigskip
\bigskip


\vfill

{\narrower
 \small
 \centerline{\textsc{Abstract}}
 \medskip\noindent
 Let $X$ be a smooth projective variety
 over a finitely generated field $K$ of characteristic~$0$
 and fix an embedding $K \subset \CC$.
 The Mumford--Tate conjecture is a precise way of saying that
 certain extra structure on the $\ell$-adic \'etale cohomology groups of~$X$
 (namely, a Galois representation)
 and
 certain extra structure on the singular cohomology groups of~$X$
 (namely, a Hodge structure)
 convey the same information.

 The main result of this paper says that if $A_1$ and~$A_2$ are
 abelian varieties (or abelian motives) over~$K$,
 and the Mumford--Tate conjecture holds for
 both~$A_1$ and~$A_2$, then it holds for $A_1 \times A_2$.
 These results do not depend on the embedding $K \subset \CC$.
 \par}

\vfill

\section{Introduction} 

\paragraph{} 
Let $A$~be an abelian variety over a finitely generated field $K \subset \CC$.
Denote with $\bar K$ the algebraic closure of $K$ in~$\CC$.
If $\ell$ is a prime number,
we write $\Hl^1(A)$ for the $\ell$-adic cohomology group
$\HH_{\et}^1(A_{\bar K}, \QQl)$.
Similarly, we write $\HB^1(A)$
for the singular cohomology group $\HH_{\sing}^1(A(\CC), \QQ)$.
There is a natural isomorphism $\Hl^1(A) \cong \HB^1(A) \otimes \QQl$
of vector spaces.

The vector space $\Hl^1(A)$ carries a Galois representation
$\rho_\ell \colon \Gal(\bar K/K) \to \GL(\Hl^1(A))$,
while $\HB^1(A)$ carries a Hodge structure.
This Hodge strucutre may be described
by a representation $\rho \colon \GB(A) \to \GL(\HB^1(A))$,
where $\GB(A)$ is the \emph{Mumford--Tate group} of~$A$
(see \cref{HS}).

Write $\Gl(A)$ for the Zariski closure of the image of~$\rho_\ell$,
and $\Glc(A)$ for the identity component of~$\Gl(A)$.
The \emph{Mumford--Tate conjecture} expresses the expectation that
the comparison isomorphism $\Hl^1(A) \cong \HB^1(A) \otimes \QQl$
identifies $\Glc(A)$ with~$\GB(A) \otimes \QQl$.
This conjecture is still wide open.
\looseness=-1

\paragraph{Main theorem} 
The goal of this article is \cref{mtcaxa}:

\smallskip

{\narrower\it\noindent
 Let $A_1$ and~$A_2$ be two abelian varieties
 over a finitely generated field $K \subset \CC$.
 If the Mumford--Tate conjecture is true for $A_1$ and~$A_2$,
 then it is true for $A_1 \times A_2$.
 \par}

\medskip

\nobreak
\noindent
In fact, in \cref{mtc-abelian-motives-tannakian-subcategory}
we prove the more general statement that
the full subcategory of abelian motives over~$K$
consisting of motives for which the
Mumford--Tate conjecture holds
is a subcategory that is closed under
direct sums, tensor products, duals, and taking direct summands.

\newpage

\begin{remark} 
 \begin{enumerate}
  \item Observe that the conclusion of the theorem
   is not a formal consequence of the assumption:
   Suppose that $G'$ is a group, with two representations
   $\rho_1 \colon G' \to \GL(V_1)$
   and
   $\rho_2 \colon G' \to \GL(V_2)$.
   Let $G_1$ (resp.~$G_2$) be the image of~$\rho_1$ (resp.~$\rho_2$).
   Write $\rho$ for $\rho_1 \oplus \rho_2$,
   and let $G$ be the image of~$\rho$.
   Then $G$ is a subgroup of $G_1 \times G_2$,
   and the projection of~$G$ onto $G_1$ (and~$G_2$)
   is surjective.
   However, $G \subset G_1 \times G_2$
   may be anything, ranging from the diagonal
   (\emph{e.g.}, if $V_1 \cong V_2$)
   to the full product
   (\emph{e.g.}, if $G_1 \not\cong G_2$ and both groups are simple).

   In the context of the main theorem we have
   \[
    \Glc(A_1 \times A_2) \subset \Glc(A_1) \times \Glc(A_2)
    \cong (\GB(A_1) \times \GB(A_2)) \otimes \QQl
    \supset \GB(A_1 \times A_2) \otimes \QQl,
   \]
   and there is no \emph{a priori} formal reason why
   $\Glc(A_1 \times A_2)$ and $\GB(A_1 \times A_2) \otimes \QQl$
   should be the same subgroup.
  \item \label[remark]{Goursat}
   The situation above is exactly the setup where Goursat's lemma applies:
   we have two groups~$G_1$ and~$G_2$
   and a subgroup $G' \subset G_1 \times G_2$
   such that the projections $\pi_i \colon G \to G_i$
   are surjective ($i = 1,2$).
   Let $N_1$ be the kernel of~$\pi_2$,
   and $N_2$ the kernel of $\pi_1$.
   Goursat's lemma is the observation that
   one may identify $N_i$ with a normal subgroup of~$G_i$,
   and the image of~$G'$ in $G_1/N_1 \times G_2/N_2$
   is the graph of an isomorphism $G_1/N_1 \to G_2/N_2$.

   This lemma is also true in the context of algebraic groups;
   a fact that we will need later on.
   We leave the proofs of these statements as an exercixe to the reader.
 \end{enumerate}
\end{remark}

\begin{remark}
 Vasiu~\cite{Va08} proves a similar result to \cref{mtcaxa}
 although he has to exclude the case where $A_1$ or~$A_2$
 has a Mumford--Tate group with a simple factor of type~$\DtD_4^\HQ$.
 His proof is long and very technical,
 and I do not claim to fully grasp the details.
 His global strategy is similar to the one employed below;
 and the reason that we can now prove the stronger claim
 is mostly due to the results of~\cite{Co17} (building on~\cite{Kisin_modp}).
\end{remark}

\paragraph{Strategy of the proof} 
\label{strategy}
\begin{enumerate}
 \item As a first step, we linearise the category of abelian varieties
  into so called \emph{abelian motives} (in the sense of Andr\'e~\cite{An95},
  or motives for absolute Hodge cycles).
  We obtain a semisimple Tannakian category,
  allowing us to apply the toolkit
  of representation theory of reductive linear groups.
 \item From work of several people (notably Piatetski-Shapiro,
  Deligne, Andr\'e, and Faltings) we know that for any abelian motive~$M$
  the group $\Glc(M)$ is reductive,
  and we have an inclusion $\Glc(M) \subset \GB(M) \otimes \QQl$.
 \item We then prove that the connected component of the centre of~$\Glc(A)$ is
  isomorphic to the connected component of the centre of~$\GB(A) \otimes \QQl$.
  For this we employ \emph{\cm~motives}, and
  reduce the claim to the Mumford--Tate conjecture for \cm~abelian varieties,
  which is known by work of Pohlmann~\cite{Pohl68}.
 \item The next step consists of replacing the abelian variety~$A_i$ ($i = 1,2$)
  by the motive~$M_i$ that corresponds---via the Tannakian formalism---with
  the adjoint representation of $\GB(A_i)^\ad$.
  It suffices to prove the Mumford--Tate conjecture for $M_1 \oplus M_2$.
 \item By general considerations,
  we may assume that $M_1$ and~$M_2$ are irreducible motives.
  In particular, the Mumford--Tate groups~$\GB(M_1)$ and~$\GB(M_2)$
  are $\QQ$-simple adjoint groups.
  In addition, we assume that
  $\Glc(M_1 \oplus M_2) \subsetneq \Glc(M_1) \times \Glc(M_2)$.
 \item We use Goursat's lemma (see \cref{Goursat}) and
  results from~\cite{Co17} to show that
  for all prime numbers~$\ell$ we have $\Hl(M_1) \cong \Hl(M_2)$.
  From this we deduce that there is a canonical isomorphism
  $\End(M_1) \cong \End(M_2)$.
 \item The remainder of the proof consists of applying
  a construction of Deligne to~$M_1$ and~$M_2$
  that is reminiscent of the Kuga--Satake construction for K3~surfaces.
  As a result we acquire two abelian varieties~$\tilde A_1$ and~$\tilde A_2$,
  and our job is to show that the isomorphism $\Hl(M_1) \cong \Hl(M_2)$
  lifts to an isomorphism $\Hl^1(\tilde A_1) \cong \Hl^1(\tilde A_2)$.
 \item Once that is done,
  we apply Faltings's theorem, to deduce that $\tilde A_1$ and~$\tilde A_2$
  are isogenous abelian varieties.
  This in turn implies $\GB(\tilde A_1) \cong \GB(\tilde A_2)$.
  In particular $\GB(M_1 \oplus M_2) \subset \GB(M_1) \times \GB(M_2)$
  is the diagonal,
  and therefore $\Glc(M_1 \oplus M_2) \cong \GB(M_1 \oplus M_2) \otimes \QQl$.
  Hence we win!
\end{enumerate}

\paragraph{Notation} 
For any field~$K$,
we denote with~$\Gamma_K$ the absolute Galois group~$\Gal(\bar K/K)$.

\paragraph{Acknowledgements} 
My warmest thanks go to my supervisor Ben Moonen.
Our countless discussions and his many detailed explanations and corrections
have been of immense importance for this article.
I also thank Rutger Noot for a very inspiring discussion of this subject.
This article also benefited from the extensive feedback on my PhD~thesis
that I received from Anna Cadoret, Pierre Deligne, Bas Edixhoven,
Milan Lopuha\"a, Rutger Noot, and Lenny Taelman.
I thank Netan Dogra, Carel Faber, Salvatore Floccari, Joost Nuiten,
and Frans Oort for their interest and useful comments.

\section{Hyperadjoint objects in Tannakian categories} 

\begin{readme} 
 In representation theory the adjoint representation is very important.
 Via the Tannakian formalism
 we port adjoint representations to Tannakian categories---%
 which are categories that behave like
 the category of representations of a linear algebraic group.
 (For a good overview of the Tannakian formalism, see~\cite{Breen_TanCats}.
 For details we refer to~\cite{Deligne_CatTan} and~\cite{DMOS}.)

 This leads to the definition of hyperadjoint objects in Tannakian categories.
 We study some of their properties in \cref{ha-props}.
 This section ties into the proof of the main theorem
 because we will replace abelian varieties~$A$ by the motive
 that corresponds to the adjoint representation of the
 motivic Galois group of~$A$.
 (See also the strategy in \cref{strategy}.)
\end{readme}

\paragraph{} 
Let $Q$ be a field of characteristic~$0$, and
let $T$ be a $Q$-linear symmetric monoidal category.
Let $R$ be a $Q$-algebra, and
denote with $\Proj_R$ the category of
finitely generated projective $R$-modules.
An \emph{$R$-valued fibre functor} of~$T$
is a $Q$-linear monoidal functor $T \to \Proj_R$
that is faithful and exact.
We denote
the groupoid of fibre functors $T \to \Proj_R$
with $\Fib(T)_R$.

\paragraph{} 
Let $Q$ be a field of characteristic~$0$.
A Tannakian category over~$Q$
is a $Q$-linear rigid abelian symmetric monoidal category
with an isomorphism $Q \stackrel\sim\to \End(1)$ such that
for every object $V \in T$ the following equivalent conditions hold:
\begin{enumerate*}[label=(\textit{\roman*})]
 \item there exists an integer~$n$ such that $\bigwedge^n V = 0$; or
 \item $\dim(V)$ is an integer.
\end{enumerate*}
(See \S1.2 and th\'eor\`eme~7.1 of~\cite{Deligne_CatTan}.)
The exterior power $\bigwedge^n V$ is defined in the usual way
in terms of $\bigotimes^n V$ and antisymmetrisation.
The dimension of~$V$ is defined as the trace of the identity morphism on~$V$,
in other words, $\dim(V)$ is the composition of the natural morphisms
$\delta$ (unit) and $\ev$ (counit):
$1 \stackrel\delta\longto V^\star \otimes V \stackrel\ev\longto 1$.
Th\'eor\`eme~7.1 of~\cite{Deligne_CatTan}
shows that the two conditions listed above are equivalent to the existence of
a $Q$-algebra~$R$ and a fibre functor $T \to \Proj_R$.

\paragraph{} 
Let $T$ be a Tannakian category over a field~$Q$ of characteristic~$0$.
For a $Q$-algebra~$R$, recall that $\Fib(T)_R$ is the groupoid of
fibre functors $T \to \Proj_R$.
It turns out that $\Fib(T)$ is an algebraic stack over~$Q$.
In fact, if $\alpha \colon Q \to R$ is a $Q$-algebra,
and $\omega \colon T \to \Proj_R$ is a fibre functor,
then the stack $\alpha^*\Fib(T)$ is isomorphic to $\BB G = [\Spec(R)/G]$,
where $G$ is the affine group scheme $\iAut^\otimes(\omega)$ over~$R$.
This observation
(together with the fact that such fibre functors exist)
makes $\Fib(T)$ into a gerbe.

A representation of~$\Fib(T)$ is a cartesian functor $\Fib(T) \to \Proj$,
in other words, a collection of functors $\Fib(T)_R \to \Proj_R$
that is functorial in~$R$.
The category of representations of~$\Fib(T)$ is denoted $\Rep(\Fib(T))$,
and the evaluation functor $T \to \Rep(\Fib(T))$,
given by $V \mapsto (\omega \mapsto \omega(V))$ is an equivalence.
This is one half of the statement of Tannaka duality.
The other half is the converse statement:
if $G$ is an affine gerbe over~$Q$,
then $G$ is naturally isomorphic to $\Fib(\Rep(G))$.

\begin{definition} 
 Let $T$ be a Tannakian category over a field~$Q$ of characteristic~$0$.
 Assume that $T$ is finitely generated (hence generated by one object).
 The \emph{adjoint object} in~$T$ is the object
 (well-defined up to isomorphism)
 that corresponds with the collection of functors
 $\Fib(T)_R \to \Proj_R$ given by $\omega \mapsto \Lie(\iAut^\otimes(\omega))$
 via the Tannakian formalism described above.
 
 Note: Since $T$ is finitely generated,
 the group scheme $\iAut^\otimes(\omega)$ is of finite type,
 and therefore $\Lie(\iAut^\otimes(\omega))$ is finitely generated.
\end{definition}

\begin{notation} 
 Let $T$ be a Tannakian category over a field~$Q$ of characteristic~$0$.
 If $V$ is an object of~$T$,
 then $V^\ad$ denotes the adjoint object
 of the Tannakian subcategory $\Tangen{V} \subset T$ generated by~$V$.
\end{notation}

\paragraph{} 
Let $T$ be a Tannakian category over a field~$Q$ of characteristic~$0$.
If $V$ is an object in~$T$,
inductively define a sequence of objects by $V^{(0)} = V$,
and $V^{(i+1)} = (V^{(i)})^\ad$ for $i \in \ZZ_{\ge0}$.
Observe that for $i \ge 1$ the object~$V^{(i+1)}$ is a quotient of~$V^{(i)}$,
and therefore $\dim V^{(i+1)} \le \dim V^{(i)}$.
Since $V$ is finite-dimensional
this sequence stabilises at an object $V^{(\infty)}$.

\begin{definition} 
 \label{ha-obj}
 Retain the notation of the preceding paragraph.
 We call the object $V^{(\infty)}$ the \emph{hyperadjoint object}
 associated with~$V$, and we denote it with~$V^\ha$.
 We say that an object $V \in T$ is \emph{hyperadjoint} if $V \cong V^\ha$
 (or equivalently, if $V \cong V^\ad$).
\end{definition}

\begin{remark} 
 Let $T$ be a Tannakian category over a field~$Q$ of characteristic~$0$.
 The constructions $V \rightsquigarrow V^\ad$ and $V \rightsquigarrow V^\ha$
 are not functorial.
 They do not in general commute with
 tensor functors between Tannakian categories.
 Also, the constructions are not in general compatible with direct sums.
 Note that the definitions are such that
 if $V \ne 0$ is a hyperadjoint object in~$T$
 then $V \oplus V$ is not hyperadjoint.
 Caveat emptor!

 On a more positive note, the following remark explains that in this paper
 these constructions are very manageable.
 \Cref{ha-props} also lists some natural properties of these constructions.
\end{remark}

\begin{remark} 
 In this paper we always have $V^\ha = V^{(2)}$
 for all objects that are of interest to us.
 The reason for this is that all the objects we encounter live in
 Tannakian (sub)categories that are semisimple,
 and therefore the associated groups (or gerbes) are reductive.
 Now suppose that $G$ is a reductive group,
 with a faithful representation $V \in \Rep(G)$.
 After the first step, we have the object $V^{(1)} = V^\ad = \Lie(G)$.
 Since $G$ is reductive, we have a short exact sequence
 $0 \to \Zentrum(G) \to G \to G^\ad \to 0$,
 and $\Lie(G) = \Lie(\Zentrum(G)) \oplus \Lie(G^\ad)$.
 Observe that $\Lie(\Zentrum(G))$ is isomorphic
 to a number of copies of the trivial representation of~$G$,
 and therefore $G^\ad$ is the group associated with~$V^{(1)}$.
 We conclude that $V^{(2)} = \Lie(G^\ad)$,
 which is a faithful representation of~$G^\ad$,
 and therefore $V^\ha = V^{(2)}$.
\end{remark}

\begin{remark} 
 I do not know of an intrinsic way to define
 adjoint and hyperadjoint objects in a finitely generated Tannakian category.
 Given the universal nature of the adjoint representation,
 I expect that it is possible to give a definition
 without using the Tannakian formalism to pass to algebraic groups or gerbes.
 Such a definition might also lead to intrinsic proofs of several properties,
 such as those in the following \namecref{ha-props}.
\end{remark}

\begin{lemma} 
 \label{ha-props}
 Let $T$ be a Tannakian category over a field~$Q$ of characteristic~$0$.
 Let $V$ be an object of~$T$, and $W$ an object of~$\Tangen{V}$.
 \begin{enumerate}
  \item We have $W^\ad \in \Tangen{V^\ad}$, and $W^\ha \in \Tangen{V^\ha}$.
  \item If $V$ is a direct sum of hyperadjoint objects,
   then $\Tangen{V}$ is semisimple.
  \item If in addition $W$ is hyperadjoint, then $W$ is a direct summand of~$V$.
 \end{enumerate}
  Suppose that $V = V_1 \oplus V_2$, with $V_1,V_2 \in T$.
 \begin{enumerate}[resume]
  \item Then $V^\ad$ is a subobject of $V_1^\ad \oplus V_2^\ad$,
   and in particular an object of $\Tangen{V_1^\ad \oplus V_2^\ad}$.
  \item For all $i \in \ZZ_{\ge0}$,
   we have $V^{(i+1)} \in \Tangen{(V_1^{(i)} \oplus V_2^{(i)})^\ad}
   \subset \Tangen{V_1^{(i+1)} \oplus V_2^{(i+1)}}$.
  \item The object $V^\ha$ is a direct summand of $V_1^\ha \oplus V_2^\ha$.
   \label[lemma]{ha-summands}
 \end{enumerate}
 \begin{proof}
  We explain the proof under the assumption
  that there is a fibre functor $\omega \in \Fib(T)_Q$.
  Let $G$ be the group $\iAut^\otimes(\omega|_{\Tangen{V}})$
  and denote with $H$ the group $\iAut^\otimes(\omega|_{\Tangen{W}})$.
  By assumption there is a surjective map $G \onto H$.
  \begin{enumerate}
   \item Since $G$ and~$H$ are groups over the field~$Q$ of characteristic~$0$,
    the map $G \onto H$ induces a surjection $\Lie(G) \onto \Lie(H)$.
    This proves the first claim; the second follows by induction.
   \item If $V$ is hyperadjoint, then $G$ is semisimple and thus reductive.
    The general case---where $V$ is a direct sum of hyperadjoint objects---%
    follows from Goursat's lemma in the context of algebraic groups
    (see \cref{Goursat}).
   \item Both $G$ and $H$ are adjoint semisimple, and $H$ is a quotient of~$G$.
    Thus $H$ is a factor of~$G$.
  \end{enumerate}
  Let $G_i$ be the group $\iAut^\otimes(\omega|_{\Tangen{V_1}})$.
  There is a natural map $G \into G_1 \times G_2$,
  and its composition with the projection onto~$G_1$ or~$G_2$ is surjective.
  \begin{enumerate}[resume]
   \item There is a natural map $\Lie(G) \into \Lie(G_1) \oplus \Lie(G_2)$.
   \item Inductively apply point~1 and the preceding point.
   \item Apply the preceding point and point~3.
  \end{enumerate}
 \end{proof}
\end{lemma}

\section{Fractional Hodge structures} \label{HS} 

\begin{readme} 
 Following \cite{Del_ShimVar} we use the notion of fractional Hodge structures.
 We also define the slightly more general notion of a
 fractional pre-Hodge structure.
 This section does not contain anything original,
 but only introduces these concepts because they will prove useful
 in understanding Deligne's construction (\cref{delignes-construction}).
\end{readme}

\begin{definition} 
 Let $R \subset \RR$ be a subring (typically $\ZZ$,~$\QQ$, or~$\RR$).
 A \emph{fractional pre-Hodge structure} over~$R$
 consists of a free $R$-module~$V$ of finite rank,
 and a decomposition
 $
  V \otimes \CC \cong \bigoplus_{p,q \in \QQ} V^{p,q}
$
 over~$\CC$, such that $V^{p,q} = \overline{V^{q,p}}$.

 We denote the category of fractional pre-Hodge structures over~$R$
 with $\FpHS_R$.
\end{definition}

\paragraph{} 
Let $V$ be a fractional pre-Hodge structure over a ring $R \subset \RR$.
For $p,q \in \QQ$, we denote with $h^{p,q}(V)$ the dimension of $V^{p,q}$.
We say that $V$ is \emph{pure} of \emph{weight}~$n \in \QQ$
if $h^{p,q}(V) \ne 0 \implies p + q = n$.
A \emph{fractional Hodge structure} is a fractional pre-Hodge structure
that is the direct sum of pure fractional pre-Hodge structures.
A \emph{pre-Hodge structure}~$V$ (without the adjective \emph{fractional})
is a fractional pre-Hodge structure
for which $h^{p,q}(V) \ne 0 \implies p,q \in \ZZ$.
If $V$ is both a fractional Hodge structure and a pre-Hodge structure,
then $V$ is a \emph{Hodge structure}, in the classical sense of the word.

\paragraph{} 
Let $\DelS$ denote the Deligne torus~$\Res_{\CC/\RR} \Gm$.
Recall that a Hodge structure over~$R$ is completely described
by a representation $h \colon \DelS \to \GL(V)_\RR$, as follows:
for $z \in \DelS(\CC)$ and $v \in V^{p,q}$
we put $h(z) \cdot v = z^{-p}\bar z^{-q}v$.
Composing $h$ with the map $x \mapsto x^k \colon \DelS \to \DelS$
amounts to relabeling $V^{p,q}$ as~$V^{kp,kq}$.

Put $\tilde \DelS = \lim_\NN \DelS$,
where $\NN$ is ordered by divisibility,
and for $m \divides n$ we take the transition map $\DelS \to \DelS$
given by $x \mapsto x^{n/m}$.
Then $\tilde\DelS$ is a pro-algebraic group scheme,
and the category of fractional pre-Hodge structures over~$\RR$
is equivalent to $\Rep(\tilde\DelS)$.

\begin{definition} 
 Let $V$ be a fractional pre-Hodge structure over a ring $R \subset \RR$.
 The \emph{Mumford--Tate group} of~$V$
 is the smallest algebraic subgroup $\GB(V) \subset \GL(V)$ over~$R$
 such that $h \colon \tilde\DelS \to \GL(V)_\RR$
 factors through $\GB(V)_\RR \subset \GL(V)_\RR$.

 Alternatively, let $\omega \colon \FpHS_R \to \Proj_R$
 be the forgetful functor.
 Then $\omega$ is a fibre functor,
 and $\GB(V) = \iAut^\otimes(\omega|_{\Tangen{V}})$.
\end{definition}

\paragraph{} 
\label{hodge-cocharacter}
Let $V$ be a pre-Hodge structure over a ring $R \subset \RR$.
Recall that this pre-Hodge structure is described by a morphism
$h \colon \DelS \to \GB(V)_\RR$.
Denote with $\mu_0$ the cocharacter $\Gm[\CC] \to \DelS_\CC$
given by $z \mapsto (z,1)$ on $\CC$-valued points.
The composite morphism $\mu_h = h_\CC \circ \mu_0$
is called the \emph{Hodge cocharacter} of~$V$.
If there is no confusion, we will write $\mu$ for~$\mu_h$.

\begin{lemma} 
 Let $V$ be a fractional pre-Hodge structure over~$\QQ$.
 The Mumford--Tate group of~$V$ is a torus
 if and only if
 $V$ is a free module of rank~$1$ over a
 commutative semisimple algebra $E \subset \End_{\FpHS_\QQ}(V)$.
 \begin{proof}
  Assume that $\GB(V)$ is a torus.
  Let $T \subset \GL(V)$ be a maximal torus containing $\GB(V)$.
  Then $E = \End_{\Vect_\QQ}(V)^T \subset \End_{\Vect_\QQ}(V)^{\GB(V)}
  = \End_{\FpHS_\QQ}(V)$ is a commutative semisimple algebra
  and $V$ has rank~$1$ over~$E$.

  Conversely, suppose that $V$ is free of rank~$1$ over some
  commutative semisimple algebra $E \subset \End_{\FpHS_\QQ}(V)$.
  Then $\GB(V) \subset \Res_{E/\QQ}\Gm = \Res_{E/\QQ} \GL_E(V)$ is a torus.
 \end{proof}
\end{lemma}

\begin{definition} 
 \label{cm-hodge-structure}
 A fractional pre-Hodge structure~$V$ over a ring $R \subset \RR$
 is called \emph{of \cm~type}
 (or a \emph{fractional \cm~pre-Hodge structure})
 if the Mumford--Tate group $\GB(V)$ is a torus.
\end{definition}

\begin{lemma} 
 \label{cm-hs-cat}
 The full subcategory of $\FpHS_\QQ$
 consisting of (classical) \cm~Hodge structures over~$\QQ$
 is a Tannakian subcategory
 generated by Hodge structures of the form $\HB^1(A)$
 where $A$ is a complex abelian variety of \cm~type.
 \begin{proof}
  See \S2 of~\cite{Andre_remarque}.
 \end{proof}
\end{lemma}

\section{Abelian motives} 

\begin{readme} 
 To give ourselves access to
 the strength and flexibility of representation theory,
 we linearise the category of abelian varieties,
 yielding a category of so-called abelian motives.
 As a consequence of work of Deligne and Andr\'e
 this category is very tractable (see \cref{hodge-is-motivated}).
\end{readme}

\paragraph{} 
The category of (pure) motives developed by Andr\'e~\cite{An95}
is very suitable for the problems at hand.
Alternatively, we could use motives for absolute Hodge cycles;
this would not influence the statements of results or proofs.
In the next paragraph we recall the definition of Andr\'e.
\looseness=-1

\paragraph{} 
Let $K$ be a subfield of~$\CC$,
and let $X$ be a smooth projective variety over~$K$.
A class $\gamma$ in $\HB^{2i}(X)$ is called
a \emph{motivated} cycle of degree~$i$
if there exists an auxiliary smooth projective variety~$Y$ over~$K$
such that $\gamma$ is of the form $\pi_*(\alpha \cup \star\beta)$,
where $\pi \colon X \times Y \to X$ is the projection,
$\alpha$ and~$\beta$ are algebraic cycle classes in $\HB^*(X \times Y)$,
and $\star\beta$ is the image of~$\beta$ under the Hodge star operation.
(Alternatively, one may use the Lefschetz star operation,
see~\S1 of~\cite{An95}.)

Every algebraic cycle is motivated,
and under the Lefschetz standard conjecture the converse holds as well.
The set of motivated cycles naturally forms a graded $\QQ$-algebra.
The category of motives over~$K$, denoted~$\Mot_K$,
consists of objects $(X,p,m)$,
where $X$ is a smooth projective variety over~$K$,
$p$ is an idempotent motivated cycle on $X \times X$,
and $m$ is an integer.
A morphism $(X,p,m) \to (Y,q,n)$
is a motivated cycle~$\gamma$ of degree $n-m$ on $Y \times X$
such that $q \gamma p = \gamma$.
We denote with $\HH(X)$ the object $(X,\Delta,0)$,
where $\Delta$ is the class of the diagonal in $X \times X$.
The K\"unneth projectors $\pi_i$ are motivated cycles,
and we denote with $\HH^i(X)$ the object $(X,\pi_i,0)$.
Observe that $\HH(X) = \bigoplus_i \HH^i(X)$.
This gives contravariant functors $\HH(\_)$ and $\HH^i(\_)$
from the category of smooth projective varieties over~$K$ to~$\Mot_K$.

\begin{theorem} 
 The category $\Mot_K$ is Tannakian over~$\QQ$,
 semisimple, graded, and polarised.
 Every classical cohomology theory of smooth projective varieties over~$K$
 factors via~$\Mot_K$.
 \begin{proof}
  See th\'eor\`eme~0.4 of~\cite{An95}.
 \end{proof}
\end{theorem}

\paragraph{} 
As the preceding theorem indicates,
the category of motives that we use is designed to have realisation functors
for all the classical cohomology theories.
For each prime number~$\ell$
we obtain an $\ell$-adic realisation functor
$\Hl \colon \Mot_K \to \Rep_{\QQl}(\Gamma_K)$
to $\ell$-adic Galois representations of~$K$,
and we obtain a Betti-realisation functor
$\HB \colon \Mot_K \to \HS_\QQ$
to the category of Hodge structures over~$\QQ$.
There is a natural isomorphism $\Hl(\_) \cong \HB(\_) \otimes \QQl$
of functors from $\Mot_K$ to $\QQl$-vector spaces.

\begin{definition} 
 Let $K$ be a subfield of~$\CC$.
 The \emph{motivic Galois group} $\GG(\Mot_K)$
 is the pro-algebraic affine group scheme $\iAut^\otimes(\HB)$ over~$\QQ$
 associated with~$\Mot_K$ via the Tannakian formalism.
 If~$M$ is a motive over~$K$, then we denote with $\GG(M)$
 the affine group scheme associated with $\Tangen{M} \subset \Mot_K$;
 it is the image of $\GG(\Mot_K)$ in~$\GL(\HB(M))$.
\end{definition}

\paragraph{} 
If $K \subset L$ is an extension of subfields of~$\CC$,
then there is a natural functor $\Mot_K \to \Mot_L$.
If $K$ and~$L$ are algebraically closed, this functor is fully faithful.
If $L$ is algebraic over~$K$, there is a short exact sequence
$1 \to \GG(\Mot_L) \to \GG(\Mot_K) \to \Gal(L/K) \to 1$.
(The Galois group $\Gal(L/K)$ corresponds via the Tannakian formalism
with the subcategory of~$\Mot_K$ generated by
the objects $\HH(\Spec(K'))$ where $K'$ is an intermediate field
$K \subset K' \subset L$.)

It is not known whether $\GG(\Mot_{\bar K})$ is connected.

\begin{notation} 
 Let $M$ be a motive over a field~$K \subset \CC$.
 We denote with $\Gl(M)$ the Zariski closure of the image of
 the Galois representation $\Gamma_K \to \GL(\Hl(M))$,
 and we write $\Glc(M)$ for the identity component of~$\Gl(M)$.
 With $\GB(M)$ we denote the Mumford--Tate group
 of the Hodge structure~$\HB(M)$.
 The realisation functors induce natural injective morphisms 
 $\Gl(M) \to \GG(M)_{\QQl}$ (via the comparison isomorphism)
 and $\GB(M) \to \GG(M)$.
\end{notation}

\begin{remark} 
 In certain situations it is expected that the natural morphisms mentioned
 above are isomorphisms.
 Let us make this more precise.
 Let $M$ be a motive over a field~$K \subset \CC$.
 Assume that $K$ is finitely generated.
 The Tate conjecture predicts that invariants of~$\Gl(M)$
 (in tensor powers of $\Hl(M)$) are algebraic,
 and in particular motivated.
 Suppose that $\Gl(M)$ is reductive, so that it is determined
 by the invariant subspace in the tensor algebra on $\Hl(M)$.
 If the invariants of $\Gl(M)$ and $\GG(M)_{\QQl}$ agree,
 then the two groups are isomorphic.
 (In fact, we do not need to assume that $\Gl(M)$ is reductive,
 by the result of~\cite{Moonen_remark}.)

 Now assume that $K = \CC$.
 The Hodge conjecture predicts that
 invariants of $\GB(M)$ in the tensor algebra on~$\HB(M)$
 are algebraic, and thus motivated.
 We already know that $\GB(M)$ is reductive,
 and hence the Hodge conjecture predicts that
 $\GB(M) \into \GG(M)$ is an isomorphism.
 See \cref{hodge-is-motivated} for an example
 of a class of motives where we know ``Hodge = motivated''.

 The Tate and Hodge conjectures are naturally complemented
 by a third conjecture: the Mumford--Tate conjecture,
 that we state below.
 We leave it to the reader to verify
 that if any two of the three conjectures hold (for all motives)
 then so does the third.
\end{remark}

\begin{conjecture} 
 \label{mtc}
 Let $M$ be a motive over a finitely generated field~$K \subset \CC$.
 Fix a prime number~$\ell$.
 The Mumford--Tate conjecture for~$M$ is the statement that
 under the comparison isomorphism $\Hl(M) \cong \HB(M) \otimes \QQl$
 we have
 \[
  \MTC_\ell(M)\colon \qquad \Glc(M) \cong \GB(M) \otimes \QQl.
 \]
 We write $\MTC(M)$ for the conjecture $\forall \ell: \MTC_\ell(M)$.
\end{conjecture}

\begin{remark} 
 \begin{enumerate}
  \item Let $M$ be a motive over a field~$K \subset \CC$.
   By definition there is a smooth projective variety~$X$
   such that $M \in \Tangen{\HH(X)}$.
   By work of Serre~\cite{Serre_letter_to_Ribet}, 
   there is a finite extension $L/K$
   such that $\Gl(\HH(X_L))$ is connected for all~$\ell$.
   Since $M \in \Tangen{\HH(X)}$,
   there is a quotient map $\Gl(\HH(X_L)) \onto \Gl(M_L)$,
   and therefore $\Gl(M_L)$ is connected for all~$\ell$.
  \item \label[remark]{ext-base-field}
   If $L/K$ is a finitely generated extension field,
   then there is an isomorphism $\Glc(M_L) \cong \Glc(M)$
   (see proposition~1.3 of~\cite{Mo15}).
   Therefore, the Mumford--Tate conjecture for~$M$
   is equivalent to the Mumford--Tate conjecture for~$M_L$.
   In particular, when trying to prove this conjecture for~$M$
   we may always assume that $\Gl(M)$ is connected for all prime numbers~$\ell$.
 \end{enumerate}
\end{remark}

\begin{definition} 
 \begin{enumerate}
  \item An \emph{Artin} motive over a field~$K \subset \CC$
   is an object in the Tannakian subcategory of~$\Mot_K$
   generated by motives of the form~$\HH(\Spec(L))$,
   where $L$ is a finite field extension of~$K$.
  \item An \emph{abelian} motive over~$K$
   is an object in the Tannakian subcategory of~$\Mot_K$
   generated by Artin motives and motives of the form~$\HH(A)$,
   where $A$ is an abelian variety over~$K$.
 \end{enumerate}
\end{definition}

\begin{remark} 
 \label{ignore-Artin-motives}
 In practice, we can ignore Artin motives in this paper.
 The Mumford--Tate conjecture is trivially true for them:
 if~$M$ is an Artin motive, then both $\GB(M)$ and~$\Glc(M)$ are trivial.
 (Note that $\Gl(M)$ can be a non-trivial finite group.)

 If $M$ is an arbitrary abelian motive,
 then there is always a finite field extension~$L/K$
 and an abelian variety~$A/L$, such that $M_L \in \Tangen{A} \subset \Mot_L$.
 By \cref{ext-base-field} we know that $\MTC(M) \iff \MTC(M_L)$
 and therefore we may restrict our attention to
 motives in the Tannakian subcategory generated by abelian varieties over~$K$.
\end{remark}

\begin{theorem} 
 \label{hodge-is-motivated}
	Consider the category of motives over~$\CC$.
 The restriction of the Hodge realisation functor $\HB(\_)$
 to the subcategory of abelian motives is a full functor.
 \begin{proof}
  See th\'eor\`eme~0.6.2 of~\cite{An95}.
 \end{proof}
\end{theorem}

\begin{remark} 
 \begin{enumerate}
  \item \label[remark]{GB-id-comp}
   Let $M$ be an abelian motive over a field $K \subset \CC$.
   A corollary of \cref{hodge-is-motivated}
   is that $\GB(M)$ is the identity component of~$\GG(M)$.
   In particular, we obtain the inclusion $\Glc(M) \subset \GB(M) \otimes \QQl$.
  \item \label[remark]{GlcM-reductive}
   Let $A$ be an abelian variety over~$K$, and
   assume that $K$ is finitely generated.
   Then $\Glc(A)$ is a reductive group
   by Satz~3 in~\S5 of~\cite{Fal83}.
   (See also~\cite{Fal84}, for the case where $K$ is not a number field.)
   Suppose that $M \in \Tangen{\HH^1(A)}$.
   Then there is a surjection $\Glc(A) \onto \Glc(M)$.
   Hence $\Glc(M)$ is reductive for every abelian motive~$M$.
  \item \label[remark]{abvar-Zlc-in-ZB}
   We also know that
   $\End_{\Glc(A)}(\Hl^1(A)) \cong \End_{\GB(A)}(\HB^1(A))_{\QQl}$,
   by Satz~4 in~\S5 of~\cite{Fal83} (see also~\cite{Fal84}).
   Now observe that $\Glc(A)$ embeds into $\End_{\Vect}(\Hl^1(A))$,
   and by definition of the centre we have
   $\Zentrum(\Glc(A)) = \Glc(A) \cap \End_{\Glc(A)}(\Hl^1(A))$,
   where the intersection takes place inside $\End_{\Vect}(\Hl^1(A))$.
   Similarly we have $\Zentrum(\GB(A))_{\QQl} =
   \GB(A)_{\QQl} \cap \End_{\GB(A)}(\HB^1(A))_{\QQl}$.
   With the preceding two points in mind,
   we conclude that $\Zentrum(\Glc(A)) \subset \Zentrum(\GB(A))_{\QQl}$.
 \end{enumerate}
\end{remark}

\begin{lemma} 
 \label{abelian-motive-centre-inclusion}
 Let $M$ be an abelian motive over
 a finitely generated field $K \subset \CC$.
 Fix a prime number~$\ell$.
 Under the comparison isomorphism $\Hl(M) = \HB(M) \otimes \QQl$
 we have an inclusion of centres
 \[
  \Zentrum(\Glc(M)) \subset \Zentrum(\GB(M))_{\QQl}.
 \]
 \begin{proof}
  After replacing $K$ with a finite field extension~$L/K$
  there is an abelian variety~$A/L$
  such that $M_L \in \Tangen{\HH^1(A)}$ (see \cref{ignore-Artin-motives}).
  By \cref{ext-base-field}, we have $\Glc(M_L) = \Glc(M)$,
  and hence we may and do assume $L = K$.
  Now we have a commutative diagram
  \[
   \begin{tikzcd}
    \Zentrum(\Glc(A)) \ar[r,hook] \ar[d,hook,"a"] &
    \Glc(A) \ar[r,two heads] \ar[d,hook,"b"] &
    \Glc(M) \ar[d,hook,"c"] \\
    \Zentrum(\GB(A))_{\QQl} \ar[r,hook] &
    \GB(A)_{\QQl} \ar[r,two heads] &
    \GB(M)_{\QQl}
   \end{tikzcd}
  \]
  where $b$ and~$c$ exist by \cref{hodge-is-motivated}
  (see \cref{GB-id-comp}),
  and $a$ exists by \cref{abvar-Zlc-in-ZB}.
  (Note: we do not claim that the rows are exact.)

  Since $\Glc(A)$ and~$\Glc(M)$ are reductive (see~\cref{GlcM-reductive}),
  we find that the image of $\Zentrum(\Glc(A))$
  under the quotient map to~$\Glc(M)$
  is exactly~$\Zentrum(\Glc(M))$.
  Similarly,
  the image of~$\Zentrum(\GB(A))_{\QQl}$ in~$\GB(M)_{\QQl}$
  is~$\Zentrum(\GB(M))_{\QQl}$.
  The result follows from the commutativity of the above diagram.
 \end{proof}
\end{lemma}

\paragraph{} 
Let $M$ be a motive over a field $K \subset \CC$.
In this article,
we say that $M$ is \emph{of \cm~type} (or a \emph{\cm~motive})
if the identity component~$\GG^\circ(M)$
of the motivic Galois group~$\GG(M)$ is a torus.
The full subcategory of~$\Mot_K$ consisting of abelian \cm~motives
is a Tannakian subcategory generated by Artin motives
and motives of the form~$\HH^1(A)$,
where $A$ is a \cm~abelian variety over~$K$ (\emph{cf.}~\cref{cm-hs-cat}).

\begin{lemma} 
 \label{mtc-cm}
 Let $M$ be an abelian \cm~motive
 over a finitely generated field $K \subset \CC$.
 Then the Mumford--Tate conjecture is true for~$M$.
 \begin{proof}
  After replacing $K$ with a finite field extension~$L/K$
  there is a \cm~abelian variety~$A/L$
  such that $M_L \in \Tangen{\HH^1(A)}$ (see \cref{ignore-Artin-motives}).
  By \cite{Pohl68}, the Mumford--Tate conjecture is true for~$A$,
  and thus $\MTC(M)$ is true by \cref{ext-base-field}.
 \end{proof}
\end{lemma}

\begin{lemma} 
 \label{centre-mtc-abelian-motive}
 Let $M$ be an abelian motive over
 a finitely generated field $K \subset \CC$.
 Fix a prime number~$\ell$.
 Under the comparison isomorphism $\Hl(M) = \HB(M) \otimes \QQl$
 we have an inclusion of centres
 $\Zentrum(\Glc(M)) \subset \Zentrum(\GB(M)) \otimes \QQl$, and
 $\Zentrum(\Glc(M))^\circ = \Zentrum(\GB(M))^\circ \otimes \QQl$.
 \begin{proof}
  (This result is proven in theorem~1.3.1 of~\cite{Va08}
  and corollary~2.11 of~\cite{UY13} using different methods.)
  The first claim is \cref{abelian-motive-centre-inclusion}.
  Therefore it suffices to show that the composition
  of the natural maps
  $\Glc(M) \into \GB(M)_{\QQl} \onto \GB(M)^\ab_{\QQl}$
  is surjective.

  Let $N$ be a faithful representation of $\GB(M)^\ab$.
  By \cref{hodge-is-motivated} and the Tannakian formalism,
  we may view $N$ as a complex motive.
  Note that $N$ is an abelian \cm~motive, by construction.
  After replacing $K$ by a finitely generated extension%
  ---which is harmless by \cref{ext-base-field}---%
  we may assume that $N$ is defined over~$K$.
  Since the Mumford--Tate conjecture holds for \cm~motives (\cref{mtc-cm})
  we find $\Glc(N) = \GB(N)_{\QQl} = \GB(M)^\ab_{\QQl}$.
  To complete the proof,
  we remark that $\Glc(N)$ is exactly the image of $\Glc(M)$
  under the composite map of the preceding paragraph.
 \end{proof}
\end{lemma}

\begin{proposition} 
 \label{mtc-adjoint-motive}
 Let $M$ be an abelian motive over
 a finitely generated field $K \subset \CC$.
 Fix a prime number~$\ell$.
 Then $\MTC_\ell(M)$
 is equivalent to $\MTC_\ell(M^\ha)$,
 where $M^\ha$ is the hyperadjoint object associated with $M \in \Mot_K$,
 see \cref{ha-obj}.
 \begin{proof}
  Consider the following commutative diagram
  \[
   \begin{tikzcd}
    0 \ar[r] &
    \Zentrum(\Glc(M)) \ar[r,hook] \ar[d,hook,"a"] &
    \Glc(M) \ar[r,two heads] \ar[d,hook,"b"] &
    \Glc(M^\ha) \ar[d,hook,"c"] \ar[r] &
    0 \\
    0 \ar[r] &
    \Zentrum(\GB(M)) \otimes \QQl \ar[r,hook] &
    \GB(M) \otimes \QQl \ar[r,two heads] &
    \GB(M^\ha) \otimes \QQl \ar[r] &
    0
   \end{tikzcd}
  \]
  where the inclusions $b$ and $c$ exist by \cref{hodge-is-motivated}
  (see \cref{GB-id-comp}),
  and $a$ exists by \cref{abelian-motive-centre-inclusion}.
  Moreover, $a$ is an isomorphism on identity components,
  by \cref{centre-mtc-abelian-motive}.
  The bottom row is exact by
  the definition of hyperadjoint objects (\cref{ha-obj})
  and \cref{hodge-is-motivated}.
  
  Certainly, if $b$ is an isomorphism, then $c$ is also an isomorphism.
  Conversely, if $c$ is an isomorphism,
  then $\dim \Glc(M) = \dim \GB(M)$.
  Since $b$ is an inclusion of connected linear groups
  it must be an isomorphism.
 \end{proof}
\end{proposition}

\section{Quasi-compatible systems of Galois representations} 

\begin{readme} 
 In general it is expected that the $\ell$-adic realisations~$\Hl(M)$
 of a motive~$M$ over a finitely generated field~$K$ of characteristic~$0$
 form a compatible system of Galois representations in the sense of Serre.

 The main theorem of~\cite{Co17} shows that if $M$~is an abelian motive
 a slightly weaker result is true;
 and this is good enough for our purposes.
 The weaker condition is called quasi-compatibility
 and in this section we recall the necessary definitions and results.
\end{readme}

\paragraph{} 
Let $\kappa$ be a finite field with $q$ elements.
We denote with $F_\kappa \in \Gamma_\kappa$
the geometric Frobenius automorphism:
\emph{i.e.,} the inverse of $x \mapsto x^q$.

\begin{definition} 
 \label{frobenius-element-number-field}
 Let $K$ be a number field.
 Let $v$ be a finite place of~$K$,
 and let $K_v$ denote the completion of~$K$ at~$v$.
 Let $\bar{K}_v$ be an algebraic closure of~$K_v$, and
 let $\bar{\kappa}/\kappa$ be the extension of residue fields
 corresponding with $\bar{K}_v/K_v$.
 The inertia group, denoted~$I_v$,
 is the kernel of the natural surjection
 $\Gamma_{K_v} \onto \Gamma_\kappa$.
 The inverse image of $F_\kappa$ in $\Gamma_{K_v}$
 is called the \emph{Frobenius coset} of~$v$.
 An element $\alpha \in \Gamma_K$ is called a
 \emph{Frobenius element with respect to~$v$}
 if there exists an embedding $\bar{K} \into \bar{K}_v$
 such that $\alpha$ is the restriction
 of an element of the Frobenius coset of~$v$.
\end{definition}

\begin{definition} 
 \label{field-model}
 Let $K$ be a finitely generated field.
 A \emph{model} of~$K$ is an
 integral scheme~$X$ of finite type over~$\Spec(\ZZ)$
 together with an isomorphism between $K$ and the function field of~$X$.
\end{definition}

\begin{remark} 
 \label{model-remark}
 If $K$ is a number field
 and $R \subset K$ is an order,
 then $\Spec(R)$ is naturally a model of~$K$.
 The only model of a number field~$K$ that is normal
 and proper over~$\Spec(\ZZ)$ is $\Spec(\mathcal{O}_{K})$.
\end{remark}

\begin{definition} 
 \label{frobenius-element-model}
 Let $K$ be a finitely generated field,
 and let $X$ be a model of~$K$.
 We denote the set of closed points of~$X$ with~$X^\cl$.
 Let $x \in X^\cl$ be a closed point,
 let $K_x$ be the function field of the Henselisation of $X$ at~$x$,
 and let $\kappa(x)$ be the residue field at~$x$.
 We denote with $I_x$ the kernel of
 $\Gamma_{K_x} \onto \Gamma_{\kappa(x)}$.
 Every embedding $\bar{K} \into \bar{K}_x$
 induces an inclusion $\Gamma_{K_x} \into \Gamma_K$.
 
 Like in \cref{frobenius-element-number-field},
 the inverse image of $F_{\kappa(x)}$
 in $\Gamma_{K_x}$
 is called the Frobenius coset of~$x$.
 An element $\alpha \in \Gamma_K$ is called a
 \emph{Frobenius element with respect to~$x$}
 if there exists an embedding $\bar{K} \into \bar{K}_x$
 such that $\alpha$ is the restriction
 of an element of the Frobenius coset of~$x$.
\end{definition}

\begin{definition} 
 \label{lambda-adic-galois-representation}
 Let $K$ be a field, let $E$ be a number field
 and let $\lambda$ be a finite place of~$E$.
 A \emph{$\lambda$-adic Galois representation} of~$K$
 is a representation of $\Gamma_K$
 on a finite-dimensional $E_\lambda$-vector space~$V$
 that is continuous for the $\lambda$-adic topology.

 We denote with $\Glambda(V) \subset \GL(V)$
 the algebraic group over~$E_\lambda$ that is
 the Zariski closure of the image of $\Gamma_K$ under this representation.
 The identity component of $\Glambda(V)$ is denoted with $\Glambdac(V)$.
\end{definition}

\begin{definition} 
 \label{unramified-representation}
 Let $K$ be a field, let $E$ be a number field
 and let $\lambda$ be a finite place of~$E$.
 Let $\rho$ be a $\lambda$-adic Galois representation of~$K$.
 Let $X$ be a model of~$K$, and let $x \in X^\cl$ be a closed point.
 We say that $\rho$ is \emph{unramified at~$x$}
 if there is an embedding $\bar{K} \into \bar{K}_x$
 for which $\rho(I_x) = \{1\}$,
 where $I_x$ is the kernel of the projection
 $\Gamma_{K_x} \onto \Gamma_{\kappa(x)}$,
 as in \cref{frobenius-element-model}.
 (Remark:
 If this condition is satisfied by one embedding,
 then it is satisfied by all embeddings.)
\end{definition}

\begin{notation} 
 \label{frobenius-characteristic-polynomial}
 \gdef\Frobcharpol{\mathrm{P}}
 Let $K$ be a finitely generated field.
 Let $E$ be a number field
 and let $\lambda$ be a finite place of~$E$.
 Let $\rho$ be a $\lambda$-adic Galois representation of~$K$.
 Let $X$ be a model of~$K$,
 and let $x \in X^\cl$ be a closed point.
 Let $F_x$ be a Frobenius element with respect to~$x$.
 Assume that $\rho$ is unramified at~$x$,
 so that the element $F_{x,\rho} = \rho(F_x)$ is well-defined up to conjugation.
 For $n \in \ZZ$, we write $\Frobcharpol_{x,\rho,n}(t)$
 for the characteristic polynomial of~$F_{x,\rho}^n$.
 Note that $\Frobcharpol_{x,\rho,n}(t)$
 does not depend on the choice of~$F_x$,
 since conjugate endomorphisms have the same characteristic polynomial.
\end{notation}

\begin{definition} 
 \label{rational-galois-representation}
	Let $K$ be a finitely generated field.
 Let $E$ be a number field
 and let $\lambda$ be a finite place of~$E$.
 Let $\rho$ be a $\lambda$-adic Galois representation of~$K$.
 Let $X$ be a model of~$K$,
 and let $x \in X^\cl$ be a closed point.
 The representation~$\rho$ is said to be \emph{$E$-rational at~$x$} if
 $\rho$~is unramified at~$x$,
 and $\Frobcharpol_{x,\rho,n}(t) \in E[t]$,
 for some $n \ge 1$.
\end{definition}

\begin{definition} 
 \label{quasi-compatible-representations}
	Let $K$ be a finitely generated field.
 Let $E$ be a number field and
 let $\lambda_{1}$ and~$\lambda_{2}$ be two finite places of~$E$.
 For $i = 1,2$, let $\rho_{i}$ be
 a $\lambda_{i}$-adic Galois representation of~$K$.
 \begin{enumerate}
  \item Let $X$ be a model of~$K$,
   and let $x \in X^{\cl}$ be a closed point.
   Then $\rho_{1}$ and~$\rho_{2}$ are said to be
   \emph{quasi-compatible at~$x$}
   if $\rho_{1}$ and~$\rho_{2}$ are both $E$-rational at~$x$,
   and if there is an integer $n$ such that
   $\Frobcharpol_{x,\rho_{1},n}(t) = \Frobcharpol_{x,\rho_{2},n}(t)$
   as polynomials in $E[t]$.
  \item \label[definition]{compatwrtmodel}
   Let $X$ be a model of~$K$.
   The representations $\rho_{1}$ and~$\rho_{2}$ are
   \emph{quasi-compatible with respect to~$X$}
   if there is a non-empty open subset $U \subset X$,
   such that $\rho_{1}$ and~$\rho_{2}$ are quasi-compatible at~$x$
   for all $x \in U^{\cl}$.
  \item \label[definition]{compatallmodel}
   The representations $\rho_{1}$ and~$\rho_{2}$ are
   \emph{quasi-compatible}
   if they are quasi-compatible with respect to every model of~$K$.
 \end{enumerate}
\end{definition}

\begin{definition} 
 \label{system-galois-representations}
 Let $K$ be a field.
 With a \emph{system of Galois representations} of~$K$ we mean a triple
 $(E, \Lambda, (\rho_{\lambda})_{\lambda \in \Lambda})$,
 where
 $E$ is a number field,
 $\Lambda$ is a set of finite places of~$E$, and
 $\rho_{\lambda}$ ($\lambda \in \Lambda$)
 is a $\lambda$-adic Galois representation of~$K$.
\end{definition}

\begin{definition} 
 \label{quasi-compatible-system}
 Let $K$ be a finitely generated field.
 Let $E$ be a number field
 and let $\Lambda$ be a set of finite places of~$E$.
 Let $\rho_{\Lambda}$ be a system of Galois representations of~$K$.
 \begin{enumerate}
  \item Let $X$ be a model of~$K$.
   The system~$\rho_{\Lambda}$ is
   \emph{quasi-compatible with respect to~$X$}
   if for all $\lambda_{1}, \lambda_{2} \in \Lambda$
   the representations $\rho_{\lambda_{1}}$ and~$\rho_{\lambda_{2}}$
   are quasi-compatible with respect to~$X$.
  \item The system~$\rho_{\Lambda}$ is called
   \emph{quasi-compatible}
   if for all $\lambda_{1}, \lambda_{2} \in \Lambda$
   the representations $\rho_{\lambda_{1}}$ and~$\rho_{\lambda_{2}}$
   are quasi-compatible.
 \end{enumerate}
\end{definition}

\begin{theorem} 
 \label{abelian-motive-quasi-compatible-realisations}
 Let $K$ be a finitely generated field of characteristic~$0$.
 Let $M$ be an abelian motive over~$K$.
 Let $E$ be a subfield of~$\End(M)$,
 and let $\Lambda$ be the set of finite places of~$E$.
 Then the system $\HLambda(M)$ is a quasi-compatible system of representations.
 \begin{proof}
  See theorem~5.1 of~\cite{Co17}.
 \end{proof}
\end{theorem}

\begin{proposition} 
 \label{quasi-compatible-semisimple-isomorphic}
 Let $K$ be a finitely generated field.
 Let $E$ be a number field, and
 let $\lambda$ be a finite place of~$E$.
 Let $\rho_1$ and~$\rho_2$ be $\lambda$-adic Galois representations of~$K$.
 If $\rho_1$ and~$\rho_2$ are semisimple, quasi-compatible,
 and $\Glambda(\rho_1 \oplus \rho_2)$ is connected,
 then $\rho_1 \cong \rho_2$.
 \begin{proof}
  See proposition~6.3 of~\cite{Co17}.
 \end{proof}
\end{proposition}

\begin{proposition} 
 \label{recover-endomorphisms}
 Let $K$ be a finitely generated field.
 Let $E$ be a number field
 and let $\Lambda$ be the set of finite places of~$E$
 whose residue characteristic is different from~$\chrc(K)$.
 \def\primes{\mathscr{L}}%
 Let $\primes$ be the set of prime numbers different from~$\chrc(K)$.
 Let $\rho_\Lambda$ be a quasi-compatible system of representations of~$K$.
 Let $\rho_\primes$ be the quasi-compatible system of representations
obtained from~$\rho_\Lambda$ by restricting to $\QQ \subset E$;
 in other words, $\rho_\ell = \bigoplus_{\lambda \divides \ell} \rho_\lambda$.
 Assume that $\Gl(\rho_\ell)$ is connected for all $\ell \in \primes$.
 Fix $\lambda_0 \in \Lambda$.
 Define the field $E' \subset E$ to be the
 subfield of~$E$ generated by elements $e \in E$
 that satisfy the following condition:

 {\narrower\noindent
  There exists a model~$X$ of~$K$,
  a point $x \in X^\cl$,
  and an integer~$n \ge 1$,\\
  such that $\Frobcharpol_{x,\rho_{\lambda_0},n}(t) \in E[t]$
  and $e$ is a coefficient of $\Frobcharpol_{x,\rho_{\lambda_0},n}(t)$.
  \par}

 \noindent
 Let $\ell$ be a prime number that splits completely in~$E/\QQ$.
 If $\End_{\Gamma_K,\QQl}(\rho_\ell) \cong E \otimes \QQl$,
 then $E = E'$.
 \begin{proof}
  We restrict our attention to a finite subset of~$\Lambda$,
  namely $\Lambda_0 = \{\lambda_0\} \cup \{ \lambda \divides \ell \}$.
  Let $U \subset X$ be a non-empty open subset such that
  for all $\lambda_1, \lambda_2 \in \Lambda_0$
  the representations
  $\rho_{\lambda_1}$ and~$\rho_{\lambda_2}$
  are quasi-compatible at all $x \in U^\cl$.
  For each $x \in U^\cl$,
  let $n_x$ be an integer such that
  $\Frobcharpol_x(t) = \Frobcharpol_{x,\rho_\lambda,n_x}(t) \in E[t]$
  does not depend on $\lambda \in \Lambda_0$.

  Let $\lambda'$ be a place of $E'$ above~$\ell$.
  Let $\lambda_1$ and~$\lambda_2$ be two places
  of~$E$ that lie above $\lambda'$.
  We view $\rho_{\lambda_1}$ and~$\rho_{\lambda_2}$
  as $\lambda'$-adic Galois representation.
  Since $\ell$ splits completely in~$E/\QQ$,
  the inclusions $\QQl \subset E'_{\lambda'} \subset E_{\lambda_i}$
  are isomorphisms.
  By definition of~$E'$ we have $\Frobcharpol_x(t) \in E'[t]$.
  Therefore $\rho_{\lambda_1}$ and~$\rho_{\lambda_2}$
  are quasi-compatible $\lambda'$-adic Galois representations;
  hence they are isomorphic by \cref{quasi-compatible-semisimple-isomorphic}.
  Let $\rho_{\lambda'}$ be the $\lambda'$-adic Galois representation
  $\bigoplus_{\lambda \divides \lambda'} \rho_\lambda$.
  We conclude that
  $\End_{\Gamma_K,E'_{\lambda'}}(\rho_{\lambda'})
   \cong \Mat_{[E:E']}(E'_{\lambda'})$,
  which implies $[E:E'] = 1$.
 \end{proof}
\end{proposition}

\section{Deligne--Dynkin diagrams} 
\label{DeligneDynkin}

\begin{readme} 
 In order to streamline the discussion of Deligne's construction
 (\cref{delignes-construction}) and the proof of the main theorem
 in the final two sections,
 this section attempts to axiomatise parts of~\cite{Del_ShimVar}
 (notably \S1.3 and~\S2.3).

 We introduce Deligne--Dynkin diagrams:
 they are Dynkin diagrams with extra structure
 that capture important information about hyperadjoint abelian motives.
 The main result of this section is \cref{deldyn-local-global},
 a local-global result for
 irreducible symplectic populated Deligne--Dynkin diagrams over~$\QQ$.
\end{readme}

\paragraph{} 
Let $\Delta$ be a connected Dynkin diagram, and
let $\Delta^+$ be the extended (or affine)
Dynkin diagram associated with~$\Delta$.
Then $\Delta^+ = \Delta \sqcup \{\alpha_0\}$.
Below we depict the connected extended Dynkin diagrams,
in which $\alpha_0$ is depicted by a grey node
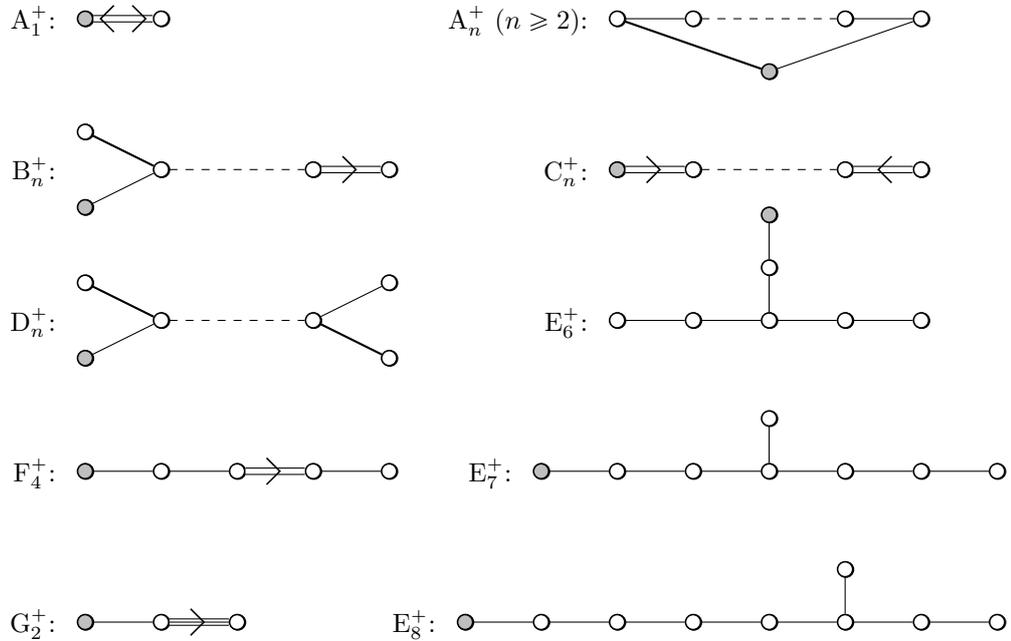
\tikz[Dynkin] \node[a0] {};.
 \[
 \begin{tikzpicture}[Dynkin] 
  \node[diagramlabel] (A1text) [label=left:{\normalsize$\DtA_1^+$:}] at (-.2,0) {};
  \node[a0] (A10) at (0,0) {};
  \node (A11) at (1,0) {};
		\draw[double distance=2pt,
		decoration={markings,
			mark=at position 0.9 with {\arrow[semithick]{Straight Barb[length=5pt]}},
			mark=at position 0.1 with {\arrowreversed[semithick]{Straight Barb[length=5pt]}}
		},
		postaction=decorate] (A10) -- (A11);

  \begin{scope}[xshift=7cm]
   \node[diagramlabel] (Antext) [label=left:{\normalsize$\DtA_n^+$ ($n \ge 2$):}] at (-.2,0) {};
   \node[a0] (An0) at (2,-0.7) {};
   \node (An1) at (0,0) {};
   \node (An2) at (1,0) {};
   \node (Annm1) at (3,0) {};
   \node (Ann) at (4,0) {};
   \draw (Annm1) -- (Ann) -- (An0);
			\draw[thick] (An0) -- (An1);
			\draw (An1) -- (An2);
   \draw[dashed] (An2) -- (Annm1);
  \end{scope}

  \begin{scope}[yshift=-2cm]
			\node[diagramlabel] (Btext) [label=left:{\normalsize$\DtB_n^+$:}] at (-.2,0) {};
   \node[a0] (B0) at (0,-.5) {};
   \node (B1) at (0,.5) {};
   \node (B2) at (1,0) {};
   \node (Bnm1) at (3,0) {};
   \node (Bn) at (4,0) {};
   \draw (B0) -- (B2);
			\draw[thick] (B2) -- (B1);
   \draw[dashed] (B2) -- (Bnm1);
   \draw[doubledynkin] (Bnm1) -- (Bn);
  \end{scope}

  \begin{scope}[xshift=7cm,yshift=-2cm]
			\node[diagramlabel] (Ctext) [label=left:{\normalsize$\DtC_n^+$:}] at (-.2,0) {};
   \node[a0] (C0) at (0,0) {};
   \node (C1) at (1,0) {};
   \node (Cnm1) at (3,0) {};
   \node (Cn) at (4,0) {};
   \draw[doubledynkin] (C0) -- (C1);
   \draw[dashed] (C1) -- (Cnm1);
   \draw[doubledynkin] (Cn) -- (Cnm1);
  \end{scope}

  \begin{scope}[yshift=-4cm]
			\node[diagramlabel] (Dtext) [label=left:{\normalsize$\DtD_n^+$:}] at (-.2,0) {};
   \node[a0] (D0) at (0,-.5) {};
   \node (D1) at (0,.5) {};
   \node (D2) at (1,0) {};
   \node (Dnm2) at (3,0) {};
   \node (Dnm1) at (4,-.5) {};
   \node (Dn) at (4,.5) {};
   \draw (D0) -- (D2);
			\draw[thick] (D2) -- (D1);
   \draw[dashed] (D2) -- (Dnm2);
   \draw[thick] (Dnm1) -- (Dnm2);
			\draw (Dnm2) -- (Dn);
		\end{scope}

  \begin{scope}[xshift=7cm,yshift=-4cm]
			\node[diagramlabel] (E6text) [label=left:{\normalsize$\DtE_6^+$:}] at (-.2,0) {};
   \node[a0] (E60) at (2,1.4) {};
   \node (E61) at (0,0) {};
   \node (E62) at (1,0) {};
   \node (E63) at (2,0) {};
   \node (E64) at (2,.7) {};
   \node (E65) at (3,0) {};
   \node (E66) at (4,0) {};
   \draw (E61) -- (E62) -- (E63) -- (E64) -- (E60);
   \draw (E63) -- (E65) -- (E66);
  \end{scope}
  \begin{scope}[xshift=6.0cm,yshift=-6cm]
   \node[diagramlabel] (E7text) [label=left:{\normalsize$\DtE_7^+$:}] at (-.2,0) {};
   \node[a0] (E70) at (0,0) {};
   \node (E71) at (1,0) {};
   \node (E72) at (2,0) {};
   \node (E73) at (3,0) {};
   \node (E74) at (3,.7) {};
   \node (E75) at (4,0) {};
   \node (E76) at (5,0) {};
   \node (E77) at (6,0) {};
   \draw (E70) -- (E71) -- (E72) -- (E73) -- (E74);
   \draw (E73) -- (E75) -- (E76) -- (E77);
  \end{scope}

  \begin{scope}[xshift=5cm,yshift=-8cm]
   \node[diagramlabel] (E8text) [label=left:{\normalsize$\DtE_8^+$:}] at (-.2,0) {};
   \node[a0] (E80) at (0,0) {};
   \node (E81) at (1,0) {};
   \node (E82) at (2,0) {};
   \node (E83) at (3,0) {};
   \node (E84) at (4,0) {};
   \node (E85) at (5,0) {};
   \node (E86) at (5,.7) {};
   \node (E87) at (6,0) {};
   \node (E88) at (7,0) {};
   \draw (E80) -- (E81) -- (E82) -- (E83) -- (E84) -- (E85) -- (E86);
   \draw (E85) -- (E87) -- (E88);
  \end{scope}

  \begin{scope}[xshift=0cm,yshift=-6cm]
   \node[diagramlabel] (F4text) [label=left:{\normalsize$\DtF_4^+$:}] at (-.2,0) {};
   \node[a0] (F40) at (0,0) {};
   \node (F41) at (1,0) {};
   \node (F42) at (2,0) {};
   \node (F43) at (3,0) {};
   \node (F44) at (4,0) {};
   \draw (F40) -- (F41) -- (F42);
   \draw[doubledynkin] (F42) -- (F43);
   \draw (F43) -- (F44);
  \end{scope}

  \begin{scope}[xshift=0cm,yshift=-8cm]
   \node[diagramlabel] (F4text) [label=left:{\normalsize$\DtG_2^+$:}] at (-.2,0) {};
   \node[a0] (G20) at (0,0) {};
   \node (G21) at (1,0) {};
   \node (G22) at (2,0) {};
   \draw[doubledynkin] (G21) -- (G22);
   \draw (G20) -- (G21) -- (G22);
  \end{scope}
 \end{tikzpicture} 
\]

\begin{definition} 
 \label{special-node}
 Let $\Delta$ be a connected Dynkin diagram, and
 let $\Delta^+ = \Delta \sqcup \{\alpha_0\}$ be the extended (or affine)
 Dynkin diagram associated with~$\Delta$.
 A node of $\Delta$ is \emph{special}
 if it is contained in the $\Aut(\Delta^+)$-orbit of~$\alpha_0$.
 See~\cref{table-deligne-dynkin-diagrams}
 for diagrams that depict which nodes are special.
\end{definition}

\begin{example} 
 \label{eg-special-nodes}
 If $\Delta$ is a connected Dynkin diagram of type~$\DtA_n$,
 then all nodes of~$\Delta$ are special.
 If $\Delta$ is of type~$\DtD_n$, then all extremal nodes are special.
 If~$\Delta$ is of type~$\DtE_8$,~$\DtF_4$, or~$\DtG_2$,
 then $\Delta$ has no special nodes.
\end{example}

\begin{definition} 
 \label{deligne-dynkin-diagram}
 A \emph{Deligne--Dynkin diagram} over a field~$Q$
 is a pair $(\Delta,\mu)$, where
 $\Delta$ is a Dynkin diagram equipped with an action of~$\Gamma_Q$,
 and $\mu$ is a subset of the special nodes
 of the connected components of~$\Delta$
 meeting each connected component in at most $1$~node.
 (Note: $\mu$~is not required to be $\Gamma_Q$-stable.)
 The Deligne--Dynkin diagram $(\Delta,\mu)$ is called \emph{irreducible}
 if $\pi_0(\Delta)$ is irreducible as $\Gamma_Q$-set,
 and $(\Delta,\mu)$ is \emph{populated} if $\mu$ meets every
 irreducible component of~$(\Delta,\mu)$.
\end{definition}

\begin{remark} 
 There is a risk of confusing terminology:
 a \emph{connected component} of~$\Delta$ or~$(\Delta,\mu)$
 will always mean a connected component of the Dynkin diagram~$\Delta$
 disregarding the $\Gamma_Q$-action;
 whereas an \emph{irreducible component} $(\Delta',\mu') \subset (\Delta,\mu)$
 is an irreducible $\Gamma_Q$-subset $\Delta'$
 of connected components of~$\Delta$ and $\mu' = \Delta' \cap \mu$.
\end{remark}

\paragraph{} 
\label{abmotdeldyn}
The reason we study Deligne--Dynkin diagrams
is that we may naturally attach such a diagram to
any hyperadjoint abelian motive~$M$ over a field $K \subset \CC$.
These assumptions on~$M$ imply that
the linear algebraic group~$\GB(M)$ over~$\QQ$ is an adjoint group.
Let $\Delta$ be the Dynkin diagram of~$\GB(M)$, and
note that $\Delta$ is naturally equipped with an action of $\Gamma_\QQ$.
Let $\mu \colon \Gm \to \GB(M)_\CC$ be the Hodge cocharacter
of the Hodge structure~$\HB(M)$
(see~\cref{hodge-cocharacter}).
This cocharacter may be identified with a subset of nodes $\mu \subset \Delta$.
Since $M$ is hyperadjoint and abelian,
the computation in \S1.2.5 of~\cite{Del_ShimVar} shows
that the nodes in~$\mu$ are special
and $(\Delta,\mu)$ is a populated Deligne--Dynkin diagram over~$\QQ$.
We call it the \emph{Deligne--Dynkin diagram of $M$}.

\paragraph{} 
\label{opposition-involution}
We recall the definition of the opposition involution on a Dynkin diagram.
Let $(R,\Phi)$ be an irreducible root system, and
let $\Delta \subset \Phi$ be a choice of positive simple roots.
Then $\Delta$ may be identified with
the vertices of the Dynkin diagram of~$(R,\Phi)$.
Let $W$ be the Weyl group of~$(R,\Phi)$, and
let $w_0$ be the longest element of the Weyl group (with respect to~$\Delta$).
Then $w_0(\Delta) = -\Delta$
and $-w_0$ defines an element~$\tau$ of $\Aut(\Delta)$:
the \emph{opposition involution}.
It is non-trivial if and only if
$\Delta$ has type $\DtA_k$ with $k \ne 1$, $\DtD_k$ with $k$~odd, or $\DtE_6$.
In these cases, $\tau$ is the unique non-trivial automorphism of~$\Delta$.
In particular, the opposition involution depends only on the type of~$\Delta$.

For non-connected Dynkin diagrams
the opposition involution is defined componentwise.

\begin{definition} 
 \label{symplectic-node}
	Let $\Delta$ be a connected Dynkin diagram with opposition involution~$\tau$,
 and let $\alpha \in \Delta$ be a special node.
 A node $\omega \in \Delta$ is called \emph{$\alpha$-symplectic}
 if $\langle \alpha, \omega + \tau(\omega) \rangle = 1$.
\end{definition}

\begin{remark} 
 The reasoning behind this terminology is best understood
 in terms of Deligne's construction
 (\cref{delignes-construction}, see also \S1.3 of~\cite{Del_ShimVar}):
 the $\alpha$-symplectic nodes of~$\Delta$ correspond precisely to
 the highest weights of certain symplectic representations.
\end{remark}

\paragraph{} 
\label{table-deligne-dynkin-diagrams}
Table~1.3.9 of~\cite{Del_ShimVar} lists
the isomorphism classes of connected Dynkin diagrams
equipped with a special node~$\alpha$.
In other words, these are
the irreducible populated Deligne--Dynkin diagrams over
an algebraically closed field~$Q = \bar Q$.

As explained in \cref{eg-special-nodes},
if the Dynkin diagram has type~$\DtD_n$, then all extremal nodes are special.
This leads to two different isomorphism classes%
---listed in the table below---%
that are labelled with $\DtD_n^\RR$ and~$\DtD_n^\HQ$.
The reader that has never seen these labels before
should not currently worry about the meaning of the superscripts~$(\_)^\RR$
and~$(\_)^\HQ$,
but the curious reader is referred to
remarque~1.3.10(ii) of~\cite{Del_ShimVar}.

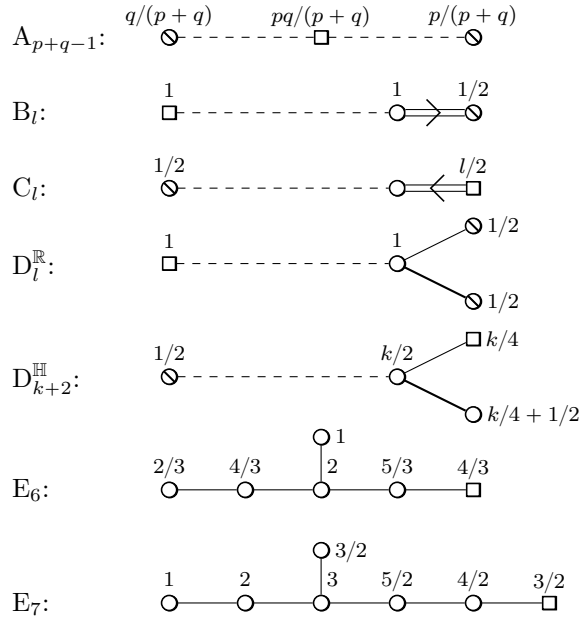
\nobreak
\noindent
\begin{minipage}{.54\textwidth} 
\[
 \begin{tikzpicture}[Dynkin,yscale=-1]
  \node[diagramlabel] (Atext) [label=right:{\normalsize$\DtA_{p+q-1}$:}] at (-1.25,0) {};
  \node[symplectic] (A1) [label={$q/(p+q)$}] at (1,0) {};
  \node[special] (Ap) [label={$pq/(p+q)$}] at (3,0) {};
  \node[symplectic] (Ak) [label={$p/(p+q)$}] at (5,0) {};
  \draw[dashed] (A1) -- (Ap);
  \draw[dashed] (Ap) -- (Ak);

  \node[diagramlabel] (Btext) [label=right:{\normalsize$\DtB_{l}$:}] at (-1.25,1) {};
  \node[special] (B1) [label={$1$}] at (1,1) {};
  \node (Blm1) [label={$1$}] at (4,1) {};
  \node[symplectic] (Bl) [label={$1/2$}] at (5,1) {};
  \draw[dashed] (B1) -- (Blm1);
  \draw[doubledynkin] (Blm1) -- (Bl);

  \node[diagramlabel] (Ctext) [label=right:{\normalsize$\DtC_{l}$:}] at (-1.25,2) {};
  \node[symplectic] (C1) [label={$1/2$}] at (1,2) {};
  \node (Clm1) at (4,2) {};
  \node[special] (Cl) [label={$l/2$}] at (5,2) {};
  \draw[dashed] (C1) -- (Clm1);
  \draw[doubledynkin] (Cl) -- (Clm1);

  \node[diagramlabel] (DRtext) [label=right:{\normalsize$\DtD_{l}^{\RR}$:}] at (-1.25,3) {};
  \node[special] (DR1) [label={$1$}] at (1,3) {};
  \node (DRlm2) [label={$1$}] at (4,3) {};
  \node[symplectic] (DRlm1) [label=right:{$1/2$}] at (5,2.5) {};
  \node[symplectic] (DRl) [label=right:{$1/2$}] at (5,3.5) {};
  \draw[dashed] (DR1) -- (DRlm2);
  \draw (DRlm1) -- (DRlm2);
  \draw[thick] (DRlm2) -- (DRl);

  \begin{scope}[yshift=0.5cm]
   \node[diagramlabel] (DHtext) [label=right:{\normalsize$\DtD_{k+2}^{\HQ}$:}] at (-1.25,4) {};
   \node[symplectic] (DH1) [label={$1/2$}] at (1,4) {};
   \node (DHlm2) [label={$k/2$}] at (4,4) {};
   \node[special] (DHlm1) [label=right:{$k/4$}] at (5,3.5) {};
   \node (DHl) [label=right:{$k/4+1/2$}] at (5,4.5) {};
   \draw[dashed] (DH1) -- (DHlm2);
   \draw (DHlm1) -- (DHlm2);
   \draw[thick] (DHlm2) -- (DHl);
  \end{scope}

  \node[diagramlabel] (E6text) [label=right:{\normalsize$\DtE_6$:}] at (-1.25,6) {};
  \node (E61) [label={$2/3$}] at (1,6) {};
  \node (E62) [label={$4/3$}] at (2,6) {};
  \node (E63) [label={\hbox{\hss\quad$2$}}] at (3,6) {};
  \node (E64) [label=right:{$1$}] at (3,5.3) {};
  \node (E65) [label={$5/3$}] at (4,6) {};
  \node[special] (E66) [label={$4/3$}] at (5,6) {};
  \draw (E61) -- (E62) -- (E63) -- (E64);
  \draw (E63) -- (E65) -- (E66);

  \begin{scope}[yshift=0.5cm]
   \node[diagramlabel] (E7text) [label=right:{\normalsize$\DtE_7$:}] at (-1.25,7) {};
   \node (E71) [label={$1$}] at (1,7) {};
   \node (E72) [label={$2$}] at (2,7) {};
   \node (E73) [label={\hbox{\hss\quad$3$}}] at (3,7) {};
   \node (E74) [label=right:{$3/2$}] at (3,6.3) {};
   \node (E75) [label={$5/2$}] at (4,7) {};
   \node (E76) [label={$4/2$}] at (5,7) {};
   \node[special] (E77) [label={$3/2$}] at (6,7) {};
   \draw (E71) -- (E72) -- (E73) -- (E74);
   \draw (E73) -- (E75) -- (E76) -- (E77);
  \end{scope}
 \end{tikzpicture}
\]
\end{minipage}\hfill
\newdimen\curparindent
\curparindent=\parindent
\begin{minipage}{.44\textwidth} 
 \parindent=\curparindent
 \divide\parindent by 3
 \multiply\parindent by 2
 \noindent
 \small
Legend (and comparison with \cite{Del_ShimVar}):
\begin{itemize}[leftmargin=12pt]
 \item[{\tikz[Dynkin] \node {};}]
  a node (bar in \cite{Del_ShimVar});
 \item[{\tikz[Dynkin] \node[special] {};}]
  the special node~$\alpha$ (circled node in \cite{Del_ShimVar});
 \item[{\tikz[Dynkin] \node[symplectic] {};}]
  an $\alpha$-symplectic node (underlined node in \cite{Del_ShimVar});
  see \cref{symplectic-node}.
\end{itemize}
The number next to a node $\omega$ is $\langle \alpha,\omega \rangle$,
where the special node~$\alpha$ and the node~$\omega$
are identified with the appropriate (co)characters.

In diagram~$\DtA_{p+q-1}$ it is the $p$-th node that is special.
The diagram~$\DtD_{k+2}^{\HQ}$ is required to satisfy $k + 2 \ge 5$.

The diagrams of type~$\DtE_6$ and~$\DtE_7$ do not have symplectic nodes.
The diagrams of type~$\DtE_8$, $\DtF_4$, and~$\DtG_2$ do not occur:
they do not have special nodes.
\end{minipage} 

\begin{definition} 
 \label{symplectic-nodes-deligne-dynkin-diagram}
 Let $(\Delta,\mu)$ be a Deligne--Dynkin diagram over a field~$Q$.
 The \emph{subset of $\mu$-symplectic nodes} attached to $(\Delta,\mu)$
 is the maximal $\Gamma_Q$-stable subset $S \subset \Delta$
 satisfying the following condition:
  For every special node $\alpha \in \mu \subset \Delta$,
  let $\Delta_\alpha$ denote the connected component of~$\Delta$
  that contains~$\alpha$.
  Then $S \cap \Delta_\alpha$ consists only of
  $\alpha$-symplectic nodes of~$\Delta_\alpha$.
  (See also \cref{table-deligne-dynkin-diagrams}.)
\end{definition}

\begin{definition} 
 A Deligne--Dynkin diagram $(\Delta,\mu)$ over~$Q$ is \emph{symplectic}
 if its subset of $\mu$-symplectic nodes has non-empty intersection
 with every irreducible component of $(\Delta,\mu)$.
\end{definition}

\begin{theorem} 
 \label{deldyn-symplectic}
 If $M$ is a hyperadjoint abelian motive over a field $K \subset \CC$,
 then the Deligne--Dynkin diagram $(\Delta,\mu)$
 associated with~$M$ is symplectic.
 \begin{proof}
  Since $(\Delta,\mu)$ depends only on~$\HB(M)$ we assume $K = \CC$.
  Because $A$ is an abelian motive,
  there is a complex abelian variety such that $M \in \Tangen{\HH^1(A)}$,
  and thus a quotient map $\GB(A) \onto \GB(M)$.
  Hence there is a map $\tilde\GG_{\mathrm{B}}(M) \to \GB(A)^{\der}$,
  where $\tilde\GG_{\mathrm{B}}(M)$ denotes
  the simply-connected cover of~$\GB(M)$.
  The representation $\GB(A) \to \GL(\HB^1(A))$ is symplectic,
  and thus its restriction to $\tilde\GG_{\mathrm{B}}(M)$ is symplectic.
  In \S1.3 and~\S2.3.7 of~\cite{Del_ShimVar}, Deligne shows that this means
  that $(\Delta,\mu)$ is symplectic.
 \end{proof}
\end{theorem}

\begin{example} 
 \label{eg-D4}
 We will now look in some detail at Deligne--Dynkin diagrams
 whose underlying Dynkin diagrams have connected components of type~$\DtD_4$.
 Let $(\Delta,\mu)$ be such a Deligne--Dynkin diagram over some field~$Q$,
 and assume that it is irreducible.
 Recall that this means that we have the following:
 \begin{itemize}
  \item a diagram~$\Delta$ consisting of a number of copies of~$\DtD_4$,
  \item an action of $\Gamma_Q$ on~$\Delta$
   that is transitive on the set~$\pi_0(\Delta)$ of connected components,
  \item a subset $\mu \subset \Delta$ of the extremal nodes
   (\emph{cf}.~\cref{eg-special-nodes}),
   such that $\mu$ meets every connected component of~$\Delta$
   in at most $1$~node.
 \end{itemize}
 Several properties of~$(\Delta,\mu)$ depend on how
 the action of $\Gamma_Q$ and the subset $\mu \subset \Delta$
 play together:
 let $\bar\mu$ denote the closure of~$\mu$
 under the action of~$\Gamma_Q$.
 Then $d = \deg_{\pi_0(\Delta)}(\bar\mu)$ takes values in $\{0,1,2,3\}$.

 The cases $d = 0$ and $d = 3$ are somewhat degenerate.
 To see this, let $S \subset \Delta$ be the subset of $\mu$-symplectic nodes.
 If $d = 0$, then $\mu = \varnothing$ and $S = \Delta$:
 in this case $(\Delta,\mu)$ is not populated.
 On the other hand, if $(\Delta,\mu)$ is populated
 then \cref{table-deligne-dynkin-diagrams} shows that
 $S$ is contained in the set~$V$ of extremal nodes of~$\Delta$.
 In fact, we have $V = \bar\mu \sqcup S$, by the maximality of~$S$.
 Hence, if $d = 3$, then $S = \varnothing$
 and $(\Delta,\mu)$ is not symplectic.
 \Cref{deldyn-symplectic,abmotdeldyn} explain why we are not interested
 in these degenerate cases.

 If $d = 1$, then $\deg(S) = 2$,
 and we say that $(\Delta,\mu)$ has type~$\DtD_4^\RR$.
 If $d = 2$, then $\deg(S) = 1$,
 and we say that $(\Delta,\mu)$ has type~$\DtD_4^\HQ$.
 We will not say more about these diagrams in this example,
 apart from mentioning that the distinction
 between diagrams of type~$\DtD_4^\RR$ and~$\DtD_4^\HQ$
 will play an important r\^ole in the proofs in the remainder of this section
 (notably the proof of \cref{deldyn-local-global}).
\end{example}

\paragraph{The type of $(\Delta,\mu)$} 
Let $(\Delta,\mu)$ be an irreducible populated Deligne--Dynkin diagram over~$Q$.
If $(\Delta,\mu)$ is symplectic,
then the type of its connected components is classical:
$\DtA_n$,~$\DtB_n$, $\DtC_n$, or~$\DtD_n$.
On the other hand, if the connected components of~$\Delta$
are of type $\DtA_n$ (resp.~$\DtB_n$, or~$\DtC_n$),
then $(\Delta,\mu)$ is symplectic.
In this case the \emph{type} of~$(\Delta,\mu)$
is $\DtA_n$ (resp.~$\DtB_n$, or~$\DtC_n$).
The case where the connected components of~$\Delta$
are of type~$\DtD_n$ requires more attention.

Assume that the connected components of~$\Delta$
are of type~$\DtD_n$, with $n \ge 5$.
For every special node $\alpha \in \mu \subset \Delta$,
let $\Delta_\alpha$ be the connected component of~$\Delta$
that contains~$\alpha$.
The pair $(\Delta_\alpha, \alpha)$ is of type~$\DtD_n^\RR$ or~$\DtD_n^\HQ$
according to the diagrams listed in \cref{table-deligne-dynkin-diagrams}.
One readily verifies that $(\Delta,\mu)$ is symplectic if and only if
one of the following conditions holds:
\begin{itemize}
 \item for every special node~$\alpha$, the pair $(\Delta_\alpha, \alpha)$
  is of type $\DtD_n^\RR$; or
 \item for every special node~$\alpha$, the pair $(\Delta_\alpha, \alpha)$
  is of type $\DtD_n^\HQ$.
\end{itemize}
We say that $(\Delta,\mu)$ is of \emph{type}~$\DtD_n^\RR$ (resp.~$\DtD_n^\HQ$)
if the former (resp.~the latter) condition holds.

Finally, if the connected components of~$\Delta$ have type~$\DtD_4$,
then the type of $(\Delta,\mu)$ was explained in \cref{eg-D4},
as follows:
Let $\bar\mu \subset \Delta$ be the $\Gamma_Q$-closure of~$\mu$
and write $d = \deg_{\pi_0(\Delta)}(\bar\mu)$.
Recall that $d \in \{1,2,3\}$,
and if $d = 3$, then $(\Delta,\mu)$ is not symplectic.
We say that $(\Delta,\mu)$ has \emph{type}~$\DtD_4^\RR$ (resp.~$\DtD_4^\HQ$)
if $d = 1$ (resp.~$d = 2$).

We conclude with two observations:
\begin{itemize}
 \item The definitions of type~$\DtD_n^\RR$ and~$\DtD_n^\HQ$
  distinguish between the cases $n = 4$ and $n \ge 5$,
  but these cases unify in the following way:
  if $(\Delta,\mu)$ is of type~$\DtD_n^\RR$ (resp.~$\DtD_n^\HQ$),
  then the subset of $\mu$-symplectic nodes of~$\Delta$
  has degree~$2$ (resp.~$1$) over~$\pi_0(\Delta)$.
 \item An irreducible symplectic populated Deligne--Dynkin diagram
  has one of the following types:
  $\DtA_n$~($n \ge 1$), $\DtB_n$ ($n \ge 2$), $\DtC_n$ ($n \ge 3$),
  $\DtD_n^\RR$ ($n \ge 4$), or~$\DtD_n^\HQ$ ($n \ge 4$);
  and all these types occur.
\end{itemize}

\begin{remark} 
 \label{restriction-remarks}
 Let $Q'/Q$ be a field extension,
 and let $(\Delta,\mu)$ be a Deligne--Dynkin diagram over~$Q$.
 By restricting the Galois action,
 one obtains a Deligne--Dynkin diagram $(\Delta,\mu)_{Q'}$ over~$Q'$.
 We make the following observations:
 \begin{enumerate}
  \item If $(\Delta,\mu)$ is irreducible or populated
   then this need not be true for $(\Delta,\mu)_{Q'}$.
   (Of course $(\Delta,\mu)_{Q'}$ will have
   irreducible components that are populated.)
  \item Related to the preceding point:
   the subset of $\mu$-symplectic nodes of $(\Delta,\mu)_{Q'}$
   may be strictly larger than
   the subset of $\mu$-symplectic nodes of $(\Delta,\mu)$.
  \item \label[remark]{restriction-remarks-type}
   If $(\Delta,\mu)$ is irreducible of type~$\DtD_4^\HQ$
   then $(\Delta,\mu)_{Q'}$ may have irreducible components of type~$\DtD_4^\RR$.
   On the other hand, if $(\Delta,\mu)$ is populated and symplectic,
   but not of type~$\DtD_4^\HQ$,
   then every irreducible component of $(\Delta,\mu)_{Q'}$ that is populated
   must have the same type as $(\Delta,\mu)$.
 \end{enumerate}
\end{remark}

\begin{lemma} 
 \label{locally-same-type}
 Let $(\Delta,\mu)$ be
 an irreducible symplectic populated Deligne--Dynkin diagram over~$\QQ$.
 Then there exists a prime number~$\ell$,
 and an irreducible component of $(\Delta,\mu)_{\QQl}$
 that has the same type as $(\Delta,\mu)$.
 \begin{proof}
  By what was said in \cref{restriction-remarks-type} this is trivial,
  unless $(\Delta,\mu)$ has type~$\DtD_4^\HQ$.
  In that case, let $\alpha \in \mu$ be a special node,
  and let $\Delta_\alpha$ be
  the connected component of~$\Delta$ that contains~$\alpha$.
  By assumption, the subset of $\mu$-symplectic nodes of $(\Delta,\mu)$
  meets $\Delta_\alpha$ in exactly one node~$s$.
  Label the remaining extremal node of~$\Delta_\alpha$ with $\beta$,
  so that $\{\alpha,\beta,s\}$ is the set of extremal nodes of~$\Delta_\alpha$.
  Note that $\bar\mu \cap \Delta_\alpha = \{\alpha,\beta\}$.
  Since $(\Delta,\mu)$ has type~$\DtD_4^\HQ$,
  there exists a special node $\alpha' \in \mu$
  and an element $g \in \Gamma_\QQ$ such that $g\alpha' = \beta$.
  By Chebotarev's density theorem
  we may assume that $g \in \Gamma_{\QQl}$ for some prime number~$\ell$.
  Let $\Delta'$ be the $\Gamma_{\QQl}$-closure of $\Delta_\alpha$,
  and take $\mu' = \Delta' \cap \mu$.
  Then $(\Delta', \mu')$ has type~$\DtD_4^\HQ$.
 \end{proof}
\end{lemma}

\paragraph{} 
An \emph{isomorphism of Deligne--Dynkin diagrams}
$\phi \colon (\Delta_1,\mu_1) \to (\Delta_2,\mu_2)$ over a field~$Q$
is a $\Gamma_Q$-equivariant isomorphism
$\phi \colon \Delta_1 \to \Delta_2$ that maps $\mu_1$ onto~$\mu_2$.
We denote the set of isomorphisms from~$(\Delta_1,\mu_1)$ to~$(\Delta_2,\mu_2)$
by $\Isom\big((\Delta_1,\mu_1),(\Delta_2,\mu_2)\big)^{\Gamma_Q}$.

Observe that there is a natural map
$\Isom\big((\Delta_1,\mu_1),(\Delta_2,\mu_2)\big)^{\Gamma_Q}
\to \Isom\big(\pi_0(\Delta_1),\pi_0(\Delta_2)\big)^{\Gamma_Q}$,
and if $f \in \Isom\big(\pi_0(\Delta_1),\pi_0(\Delta_2)\big)^{\Gamma_Q}$
then we write $\Isom_f\big((\Delta_1,\mu_1),(\Delta_2,\mu_2)\big)^{\Gamma_Q}$
for the set of
$\phi \in \Isom\big((\Delta_1,\mu_1),(\Delta_2,\mu_2)\big)^{\Gamma_Q}$
such that $\pi_0(\phi) = f$.

\begin{lemma} 
 \label{AutDeltamu}
 Let $(\Delta,\mu)$ be
 an irreducible symplectic populated Deligne--Dynkin diagram over~$Q$.
 Let $\tau$ denote the opposition involution on~$\Delta$.
 Then
 \[
  \#\Aut_\id(\Delta,\mu)^{\Gamma_Q} =
  \begin{cases}
   1 &\text{if $(\Delta,\mu)$ has type $\DtA_1$,~$\DtB_n$, $\DtC_n$,
    $\DtD_n^\HQ$, or $\DtA_n$ and $\mu$ is not fixed by~$\tau$,}\\
   2 &\text{if $(\Delta,\mu)$ has type $\DtD_n^\RR$, or $\DtA_n$ ($n \ge 2$)
    and $\mu$ is fixed by~$\tau$.}
  \end{cases}
 \]
 \begin{proof}
  If $\Delta$ has only one connected component,
  then the result may be deduced from the tables
  in~\cref{table-deligne-dynkin-diagrams}.
  If $\Delta$ has multiple components,
  then the result follows from the fact that
  a non-trivial element of $\Aut_\id(\Delta,\mu)^{\Gamma_Q}$
  must act non-trivially on every connected component of~$\Delta$.
 \end{proof}
\end{lemma}

\begin{lemma} 
 \label{locally-same-aut}
 Let $(\Delta,\mu)$ be
 an irreducible symplectic populated Deligne--Dynkin diagram over~$\QQ$.
 Then there exists
 an irreducible component $(\Delta_\lambda,\mu_\lambda)$ of $(\Delta,\mu)_{\QQl}$
 such that the natural map
 $\Aut_\id(\Delta,\mu)^{\Gamma_Q} \to
 \Aut_\id(\Delta_\lambda,\mu_\lambda)^{\Gamma_Q}$
 is an isomorphism.
 \begin{proof}
  This follows immediately from \cref{locally-same-type,AutDeltamu}.
 \end{proof}
\end{lemma}

\begin{remark} 
 \label{quadratic-remark}
 As explained in the readme at the beginning of this \namecref{DeligneDynkin},
 our aim is a local-global result for
 irreducible symplectic populated Deligne--Dynkin diagrams over~$\QQ$
 (\cref{deldyn-local-global}).
 Now recall that a quadratic extension~$F$ of a number field~$E$
 is completely determined by the set of primes of~$E$ that split in~$F$.
 This fact inspires the following notation.
\end{remark}

\paragraph{} 
\label{deg2-set}
Let $(\Delta,\mu)$ be
an irreducible symplectic populated Deligne--Dynkin diagram over~$Q$.
The action of $\Gamma_Q$ on~$\Delta$
is determined by the action on
a $\Gamma_Q$-closed subset $U(\Delta,\mu)$ of~$\Delta$
with $\deg_{\pi_0(\Delta)}(U(\Delta,\mu)) \in \{1,2\}$.
Indeed, let $S$ be the subset of $\mu$-symplectic nodes of~$(\Delta,\mu)$.
\begin{itemize}
 \item If $(\Delta,\mu)$ has type~$\DtA_n$, we take~$U(\Delta,\mu) = S$,
  which has degree~$1$ if $n = 1$ and degree~$2$ otherwise.
 \item If $(\Delta,\mu)$ has type~$\DtB_n$, we take $U(\Delta,\mu) = S$,
  which has degree~$1$.
 \item If $(\Delta,\mu)$ has type~$\DtC_n$, we take $U(\Delta,\mu) = \bar\mu$,
  which has degree~$1$.
 \item If $(\Delta,\mu)$ has type~$\DtD_n^\RR$, we take $U(\Delta,\mu) = S$,
  which has degree~$2$.
 \item If $(\Delta,\mu)$ has type~$\DtD_n^\HQ$,
  we take $U(\Delta,\mu) = \bar\mu$,
  which has degree~$1$ or~$2$.
\end{itemize}

\paragraph{} 
Let $E$ be a number field, and
let $(\Delta_1,\mu_1)$ and~$(\Delta_2,\mu_2)$ be two
irreducible symplectic populated Deligne--Dynkin diagrams over~$\QQ$
such that $\pi_0(\Delta_1)$ and $\pi_0(\Delta_2)$ are
both isomorphic to $\Hom(E,\QQbar)$ as $\Gamma_\QQ$-sets.
Fix $\Gamma_\QQ$-equivariant isomorphisms
$\pi_0(\Delta_1) \cong \Hom(E,\QQbar) \cong \pi_0(\Delta_2)$
and write $f$ for the composite map $\pi_0(\Delta_1) \to \pi_0(\Delta_2)$.
Suppose that there is an isomorphism
$\phi \in \Isom_f\big((\Delta_1,\mu_1),(\Delta_2,\mu_2)\big)^{\Gamma_\QQ}$.

We restrict the Galois action to $\Gamma_{\QQl}$,
for some prime number~$\ell$.
The irreducible components of~$(\Delta_1,\mu_1)_{\QQl}$
are in a natural way indexed by the places~$\lambda$ of~$E$
that lie above~$\ell$, so that
$(\Delta_1,\mu_1)_{\QQl} = \sqcup_{\lambda\divides\ell} (\Delta_1,\mu_1)_\lambda$.
In a similar manner, the maps~$f_\ell$ and~$\phi_\ell$
are the disjoint union of local components~$f_\lambda$ and~$\phi_\lambda$.
We will use this notation below.

\begin{proposition} 
 \label{deldyn-local-global}
 Let $E$ be a number field, and
 let $(\Delta_1,\mu_1)$ and~$(\Delta_2,\mu_2)$ be two
 irreducible symplectic populated Deligne--Dynkin diagrams over~$\QQ$
 such that $\pi_0(\Delta_1)$ and $\pi_0(\Delta_2)$ are
 both isomorphic to $\Hom(E,\QQbar)$ as $\Gamma_\QQ$-sets.
 Fix $\Gamma_\QQ$-equivariant isomorphisms
 $\pi_0(\Delta_1) \cong \Hom(E,\QQbar) \cong \pi_0(\Delta_2)$
 and write $f$ for the composite map $\pi_0(\Delta_1) \to \pi_0(\Delta_2)$.
 Suppose that for each prime number~$\ell$ there is an element
 $\psi_\ell \in
 \Isom_f\big((\Delta_1,\mu_1)_{\QQl},(\Delta_2,\mu_2)_{\QQl}\big)^{\Gamma_{\QQl}}$.
 Then there exists an isomorphism
 $\phi \in \Isom_f\big((\Delta_1,\mu_1),(\Delta_2,\mu_2)\big)^{\Gamma_\QQ}$
 and a finite place~$\lambda$ of~$E$
 such that $\phi_\lambda = \psi_\lambda$.
 \begin{proof}
  It suffices to prove the existence of~$\phi$;
  the second claim will follow automatically from \cref{locally-same-aut}.
  First of all, observe
  that $(\Delta_1,\mu_1)$ and~$(\Delta_2,\mu_2)$ have the same type,
  by \cref{restriction-remarks-type} and \cref{locally-same-type}.
  If this type is $\DtB_n$ or~$\DtC_n$,
  then the topology of the diagrams
  and the map $f \colon \pi_0(\Delta_1) \to \pi_0(\Delta_2)$
  determine a unique $\Gamma_\QQ$-equivariant map
  $\phi \in \Isom_f(\Delta_1,\Delta_2)^{\Gamma_\QQ}$.
  Note that $\phi$ maps $\mu_1$ to~$\mu_2$,
  because it does so locally.

  The other types require more bookkeeping,
  although the strategy remains the same.
  We begin by proving that
  there is an element $\phi \in \Isom_f(\Delta_1,\Delta_2)^{\Gamma_\QQ}$,
  that may or may not map $\mu_1$ to~$\mu_2$.
  By \cref{deg2-set}, a $\Gamma_\QQ$-equivariant isomorphism
  $U(\Delta_1,\mu_1) \to U(\Delta_2,\mu_2)$ lying above~$f$
  will extend to a $\Gamma_\QQ$-equivariant isomorphism
  $\Delta_1 \to \Delta_2$.
  Since the degree of~$U(\Delta_i,\mu_i)$ over~$\pi_0(\Delta_i)$
  is at most~$2$,
  the existence of global isomorphism
  $U(\Delta_1,\mu_1) \to U(\Delta_2,\mu_2)$
  will follow from Chebotarev's density theorem
  if we prove that these sets are locally isomorphic
  (\emph{cf.,} \cref{quadratic-remark}).
  \begin{itemize}
   \item If $(\Delta_1,\mu_1)$ and~$(\Delta_2,\mu_2)$ have type~$\DtA_n$,
    or~$\DtD_n^\RR$ with $n \ge 5$,
    then the topology of the diagrams
    forces that $\psi_\ell(U(\Delta_1,\mu_1))$ equals $U(\Delta_2,\mu_2)$
    for all prime numbers~$\ell$.
   \item Now assume that $(\Delta_1,\mu_1)$ and~$(\Delta_2,\mu_2)$
    have type~$\DtD_4^\RR$ or~$\DtD_4^\HQ$.
    Let $\lambda$ be a finite place of~$E$.
    There are two possibilities:
    either the extremal nodes of
    $(\Delta_1,\mu_1)_\lambda$ and $(\Delta_2,\mu_2)_\lambda$
    form $2$ orbits under the action of $\Gamma_{\QQl}$,
    or the extremal nodes form $3$ orbits.
    For both possibilities it is clear that
    $U(\Delta_1,\mu_1)_\lambda$
    is isomorphic to
    $U(\Delta_2,\mu_2)_\lambda$.
   \item Finally, assume that $(\Delta_1,\mu_1)$ and~$(\Delta_2,\mu_2)$
    have type $\DtD_n^\HQ$ with $n \ge 5$.
    Recall that the subset $U(\Delta_i,\mu_i) = \bar\mu_i$
    has degree~$1$ or~$2$.
    If these degrees are equal, for $i = 1,2$, then we are done.
    Now suppose that the degree of $\bar\mu_1$ is~$2$.
    Let $\{\alpha,\beta\} \subset \bar\mu_1$
    be two nodes that lie in the same connected component of~$\Delta_1$.
    Without loss of generality we may and do assume $\alpha \in \mu_1$.
    By Chebotarev's density theorem,
    there exists a prime number~$\ell$ and $g \in \Gamma_{\QQl}$
    such that $g\alpha = \beta$.
    By assumption, we have $\psi_\ell(\alpha) \in \mu_2$
    and $\psi_\ell(\beta) = g\psi_\ell(\alpha) \in \bar\mu_2$.
    We conclude that $\bar\mu_1$ and~$\bar\mu_2$ have the same degree.
  \end{itemize}
  We have now proven that there exists a $\Gamma_\QQ$-equivariant
  isomorphism $\phi \colon \Delta_1 \to \Delta_2$ lying over~$f$.
  It remains to prove that we can choose~$\phi$
  in such a way that it maps $\mu_1$ to~$\mu_2$.

  At the beginning of this proof, we already dealt with the cases where
  $(\Delta_1,\mu_1)$ and $(\Delta_2,\mu_2)$ have type~$\DtB_n$ or~$\DtC_n$.
  We continue with some other easy cases,
  where the topology of the diagrams forces $\phi(\mu_1) = \mu_2$.
  This happens:
  \begin{itemize}
   \item if $(\Delta_1,\mu_1)$ and $(\Delta_2,\mu_2)$
    are of type $\DtA_1$, or~$\DtD_n^\RR$ with $n \ge 5$;
   \item if the type is $\DtD_4^\RR$,
    since $\phi$ maps $U(\Delta_1,\mu_1) = S_1$
    to $U(\Delta_2,\mu_2) = S_2$ by construction;
   \item if the type is $\DtA_n$ with $n \ge 2$,
    and $\mu_1$ (and thus $\mu_2$) is fixed
    under the opposition involution.
  \end{itemize}
  The remaining cases require a bit more work:
  so far we have not needed to change our choice of~$\phi$,
  but with the remaining cases this might be necessary.

  First consider the case~$\DtD_4^\HQ$.
  Recall that $\bar\mu_i$ has degree~$2$.
  Thus there is a unique non-trivial $\Gamma_\QQ$-equivariant
  involution~$\tau_i \in \Aut_\id(\Delta_i)$
  such that $\tau_i(\bar\mu_i) = \bar\mu_i$.
  Observe that $\phi \circ \tau_1 = \tau_2 \circ \phi$.
  Let $\alpha \in \mu_1$ be a special node.
  We may assume that $\phi(\alpha) \in \mu_2$,
  by replacing $\phi$ with $\tau_2 \circ \phi$ if necessary.

  We claim that $\phi(\mu_1) = \mu_2$.
  Indeed, let $\alpha' \in \mu_1$ be a special node.
  Recall from \cref{eg-D4} that $\deg(\bar\mu_i) = 2$.
  Now there are two cases:
  either $\alpha$ and~$\alpha'$ are in the same~$\Gamma_\QQ$-orbit,
  or they are in different orbits.
  Assume that $\alpha$ and~$\alpha'$ are in the same orbit.
  The proof of this case contains the essential idea
  that will also be applied in all the remaining cases:
  There exists a $g \in \Gamma_\QQ$
  such that $g\alpha = \alpha'$.
  By Chebotarev's density theorem,
  we may assume that $g \in \Gamma_{\QQl}$
  for some prime number~$\ell$.
  Since $\phi(\alpha) \in \mu_2$, we have $\phi(\alpha) = \psi_\ell(\alpha)$,
  and we compute
  \[
   \phi(\alpha') = \phi(g\alpha) = g\phi(\alpha) = g\psi_\ell(\alpha)
   = \psi_\ell(g\alpha) = \psi_\ell(\alpha') \in \mu_2.
  \]
  Now suppose that $\alpha$ and~$\alpha'$ are in different orbits.
  Let $\Delta_\alpha$ and~$\Delta_{\alpha'}$
  be the connected components of~$\Delta$
  that contain $\alpha$ and~$\alpha'$ respectively.
  Let $s$ and~$s'$ be the $\mu$-symplectic nodes
  in~$\Delta_\alpha$ and~$\Delta_{\alpha'}$ respectively.
  Finally, let $\beta \in \Delta_\alpha$
  and~$\beta' \in \Delta_{\alpha'}$ be such that
  $\{\alpha,\beta,s\}$ is the set of extremal nodes of~$\Delta_\alpha$, and
  $\{\alpha',\beta',s'\}$ is the set of extremal nodes of~$\Delta_{\alpha'}$.
  There exists a $g \in \Gamma_\QQ$
  such that $g\alpha = \beta'$.
  By Chebotarev's density theorem,
  we may assume that $g \in \Gamma_{\QQl}$
  for some prime number~$\ell$.
  Let $\lambda \divides \ell$ be the finite place of~$E$
  that corresponds with the $\Gamma_{\QQl}$-closure of~$\Delta_\alpha$.

  By construction $\phi$ maps
  $U(\Delta_1,\mu_1) = \bar\mu_1$ to $U(\Delta_2,\mu_2) = \bar\mu_2$.
  On the other hand,
  we have $\psi_\lambda(\bar\mu_{1,\lambda}) = \bar\mu_{2,\lambda}$
  by assumption;
  thus $\phi(s) = \phi_\lambda(s)$.
  Since $\phi(\alpha) \in \mu_2$, we have $\phi(\alpha) = \psi_\lambda(\alpha)$,
  and therefore $\phi(\beta) = \psi_\lambda(\beta)$.
  We compute
  $\phi(\alpha') = \phi(g\beta) = g\phi(\beta) = g\psi_\lambda(\beta)
   = \psi_\lambda(g\beta) = \psi_\lambda(\alpha') \in \mu_2$.
  This finishes the proof in the case~$\DtD_4^\HQ$.

  Finally, suppose that the type is $\DtD_n^\HQ$ with $n \ge 5$,
  or the type is $\DtA_n$ and $\mu_i$ is not fixed by the opposition involution.
  The proof is similar to the case~$\DtD_4^\HQ$.
  There is a unique non-trivial $\Gamma_\QQ$-equivariant
  involution~$\tau_i \in \Aut_\id(\Delta_i)$.
  Once again, we have $\phi \circ \tau_1 = \tau_2 \circ \phi$.
  Let $\alpha \in \mu_1$ be a special node that is not fixed under~$\tau_1$.
  We may assume that $\phi(\alpha) \in \mu_2$,
  by replacing $\phi$ with $\tau_2 \circ \phi$ if necessary.

  As before, we claim that $\phi(\mu_1) = \mu_2$.
  Indeed, suppose that $\alpha' \in \mu_1$ is a special node, and
  let $\Delta_{\alpha'}$ be the component of~$\Delta_1$ containing~$\alpha'$.
  There exists a $g \in \Gamma_\QQ$
  such that $g\alpha \in \Delta_{\alpha'}$.
  By Chebotarev's density theorem,
  we may assume that $g \in \Gamma_{\QQl}$
  for some prime number~$\ell$.
  Now we argue as follows:
  Let $(\Delta_1,\mu_1)_\lambda$
  (resp.~$(\Delta_2,\mu_2)_\lambda$)
  be the irreducible component of~$(\Delta_1,\mu_1)_{\QQl}$
  (resp.\ $(\Delta_2,\mu_2)_{\QQl}$)
  that contains~$\alpha$ (resp.~$\phi(\alpha)$).
  Since $\phi(\alpha) \in \mu_2$, we have $\phi(\alpha) = \psi_\ell(\alpha)$
  and therefore $\phi_\lambda = \psi_\lambda$.
  This implies $\phi(\alpha') \in \mu_2$,
  by a computation similar to the one above above.
  This completes the proof of \cref{deldyn-local-global}.
 \end{proof}
\end{proposition}

\paragraph{} 
Let $M$ be a hyperadjoint abelian motive
over a finitely generated field $K \subset \CC$.
Let $(\Delta,\mu)$ be the Deligne--Dynkin diagram associated with~$M$.
Choose an embedding of $K$ into a $p$-adic field~$K_v$.
Let $C_v$ be the completion of an algebraic closure of~$K_v$,
let $C_v(i)$ denote the invertible $C_v$-module $C_v \otimes_{\QQp} \QQp(i)$,
and let $\BHT{v}$ denote the Hodge--Tate ring $\bigoplus_{i \in \ZZ} C_v(i)$.
By $p$-adic Hodge theory,
we know that there is a grading on $\Hp(M) \otimes_{\QQp} \BHT{v}$.
This determines a cocharacter~$\mu_\HT$ of
$\Gpc(M) \otimes_{\QQp} \BHT{v} \subset \GB(M) \otimes_\QQ \BHT{v}$
that is called the \emph{Hodge--Tate} cocharacter.

\begin{lemma} 
 \label{hodge-tate-cochar-deldyn}
 Retain the notation of the preceding paragraph.
 The conjugacy class of the cocharacter
 $\mu_\HT \in X_*(\GB(M) \otimes_\QQ \BHT{v})$
 and
 the conjugacy class of the Hodge cocharacter
 $\mu \in X_*(\GB(M) \otimes_\QQ \CC)$
 correspond to the same conjugacy class of
 cocharacters $X_*(\GB(M) \otimes_\QQ \QQbar)$.
 In particular, $\mu_\HT$ determines the subset $\mu \subset \Delta$.
 \begin{proof}
  By the Tannakian formalism,
  the cocharacters $\mu_\HT$ and~$\mu$ correspond to tensor functors
  \[
   \begin{aligned}
    \Hp \otimes \BHT{v} \colon \Tangen{M} &\to \grVect_{\BHT{v}} \\
    N &\mapsto \Hp(N) \otimes_{\QQp} \BHT{v}
   \end{aligned}
   \qquad\text{and}\qquad
   \begin{aligned}
    \HB \otimes \CC \colon \Tangen{M} &\to \grVect_{\CC} \\
    N &\mapsto \HB(N) \otimes_\QQ \CC.
   \end{aligned}
  \]
  There is also another tensor functor
  \begin{align*}
   \gr\HdR \colon \Tangen{M} &\to \grVect_{K} \\
   N &\mapsto \gr\HdR(N)
  \end{align*}
  that maps a motive to the associated graded
  of its algebraic de~Rham realisation.
  If $\GdR(M)/K$ denotes the algebraic group $\iAut^\otimes(\HdR)$,
  that corresponds to~$\Tangen{M}$ via the fibre functor~$\HdR$,
  then $\gr\HdR$ corresponds to a cocharacter $\mu_\dR \in X_*(\GdR(M))$.

  The groups $\GdR(M)$ and $\GG(M) \otimes_\QQ K$ are inner forms of each other.
  As a reminder:
  Recall that $\GG(M)/\QQ$ is the algebraic group $\iAut^\otimes(\HB)$,
  that corresponds to~$\Tangen{M}$ via the fibre functor~$\HB$,
  and recall from \cref{GB-id-comp} that the identity component of~$\GG(M)$
  is the Mumford--Tate group $\GB(M)$, since $M$ is an abelian motive.

  The comparison theorem between
  the Betti realisation
  and
  the algebraic de~Rham realisation
  gives a natural transformation $\HB \otimes \CC \to \gr\HdR \otimes \CC$.
  This yields an injective homomorphism of complex algebraic groups
  $\GB(M) \otimes_\QQ \CC \into \GdR(M) \otimes_K \CC$
  that maps $\mu$ to~$\mu_\dR$.

  On the other hand, $p$-adic Hodge theory compares
  the $p$-adic \'etale realisation with the algebraic de~Rham realisation,
  and gives a natural transformation
  $\Hp \otimes \BHT{v} \to \gr\HdR \otimes \BHT{v}$.
  This yields an injection
  $\GB(M) \otimes_\QQ \BHT{v} \into \GdR(M) \otimes_K \CC$
  that maps $\mu_\HT$ to~$\mu_\dR$.

  The proof concludes by an application of the following remark:
  Let $\bar Q \subset \bar Q'$ be an extension of algebraically closed fields,
  and let $G$ be an algebraic group over~$\bar Q$.
  Then the natural map $X_*(G) \to X_*(G_{\bar Q'})$ obtained by base change
  induces a bijection
  between the conjugacy classes of cocharacters of~$G$
  and the conjugacy classes of cocharacters of~$G_{\bar Q'}$.
 \end{proof}
\end{lemma}

\section{Deligne's construction} 
\label{delignes-construction}

\begin{readme} 
 Let $M$ be an irreducible hyperadjoint abelian motive over~$\CC$.
 Proposition~2.3.10 of~\cite{Del_ShimVar} provides a recipe
 to construct a complex abelian variety~$A$ (up to isogeny)
 such that $M \cong \HH^1(A)^\ha$.
 Deligne uses the language of Shimura data and Hodge theory.
 We recall this construction,
 using the terminology of abelian motives and Deligne--Dynkin diagrams.
\end{readme}

\paragraph{Preparations} 
Let $M$ be an irreducible hyperadjoint abelian motive over~$\CC$,
and let $(\Delta,\mu)$ be the Deligne--Dynkin diagram associated with~$M$.
The endomorphism algebra $E = \End(M)$ is a totally real field
(see the discussion in~\S2.3.4(a) of~\cite{Del_ShimVar}).
Note that $\Hom(E,\QQbar) \cong \pi_0(\Delta)$ as $\Gamma_\QQ$-sets.
We also recall that $(\Delta,\mu)$ is symplectic (\cref{deldyn-symplectic}).
Let $S$ be the subset of $\mu$-symplectic nodes of~$(\Delta,\mu)$.

Write $G$ for the $\QQ$-simple adjoint group $\GB(M)$,
and let $\tilde G$ be the simply-connected cover of~$G$.
For $s \in S$, let $V(s)$ be the representation of~$\tilde G_\CC$
whose highest weight corresponds to~$s$.
Since $S$ is closed under the action of~$\Gamma_\QQ$,
there exists a representation~$V$ of~$\tilde G$ over~$\QQ$
such that $V_\CC \cong \bigoplus_{s \in S} V(s)^{\oplus n}$
for suitable~$n$.

\paragraph{Choices} 
\label{choices}
Now we fix three choices:
\begin{enumerate}
 \item Choose a totally imaginary quadratic extension $F/E$.
 \item Choose a partial \cm~type~$\Phi$ for~$F$ relative to $(\Delta,\mu)$:
  a subset $\Phi \subset \Hom(F,\CC) = \Hom(F,\QQbar)$
  that maps 1-to-1 onto the complement of the image of~$\mu$
  in $\Hom(E,\QQbar) = \pi_0(\Delta)$.
 \item Choose a representation~$V$ of~$\tilde G$ as above.
\end{enumerate}

\paragraph{} 
The Hodge cocharacter $h \colon \DelS \to G_\RR$ lifts to a map
$\tilde h \colon \tilde \DelS \to \tilde G_\RR$,
endowing~$V$ with a fractional Hodge structure.
The Hodge decomposition of~$V$ may be read off from
the diagrams in~\cref{table-deligne-dynkin-diagrams}:
If $s \in S$ lies in a component of~$\Delta$ that does not meet~$\mu$,
then the type of~$V(s)$ is $\{(0,0)\}$.
If $s$ lies in a component of~$\Delta$
that contains a special node $\alpha \in \mu$,
then the type of~$V(s)$ is $\{(r,-r), (r-1, 1-r)\}$
where $r = \langle s, \alpha \rangle$ is the number
that is written next to the node~$s$
in the appropriate diagram in~\cref{table-deligne-dynkin-diagrams}.

\paragraph{} 
Let $F_S$ denote the \'etale $E$-algebra
such that $\Hom(F_S, \QQbar) \cong S$ as $\Gamma_\QQ$-sets.
Observe that the fractional Hodge structure~$V$
is canonically an $F_S$-module:
the algebra~$F_S$ acts on~$V(s)$
via the embedding $F_S \into \QQbar \subset \CC$
that corresponds with $s \in S$.

We endow $F_S$ with a fractional pre-Hodge structure:
the component $\CC^{\{s\}}$ of $F_S \otimes_\QQ \CC \cong \CC^S$
is placed in bi-degree $(0,0)$
if $s$ lies in a component of~$\Delta$ that does not meet~$\mu$;
and $\CC^{\{s\}}$ is placed in bi-degree $(1-r,r)$
if $s$ lies in a component that does meet~$\mu$,
where $r$ is the rational number from the preceding paragraph.

In a similar fashion we endow the \cm~field~$F$ with
a fractional pre-Hodge structure:
the component $\CC^\phi$ of $F \otimes \CC \cong \CC^{\Hom(F,\CC)}$
is placed in bi-degree
\[
 \begin{cases}
  (1,0) &\text{if $\phi \in \Phi$}\\
  (0,1) &\text{if $\bar\phi \in \Phi$}\\
  (0,0) &\text{otherwise.}
 \end{cases}
\]

\paragraph{} 
\label{WF}
Write $W_F$ for $F \otimes_E F_S$, and
observe that $F \otimes_E F_S$ is a fractional Hodge structure of weight~$1$.
Note that $W_F$ is of \cm~type,
since both~$F$ and $F_S$ are of \cm~type.
Put $V' = W_F \otimes_{F_S} V$.
A computation shows that $V'$
is a Hodge structure of type $\{(1,0), (0,1)\}$.
It turns out that $V'$~is a polarisable Hodge structure
(see \cite{Del_ShimVar} for details).
Thus there is a complex abelian variety~$A$ (well-defined up to isogeny)
such that $\HB^1(A) \cong V'$.
Since $W_F$ is of \cm~type,
we find that $\GB(V')^\ad = G$,
and therefore $\HH^1(A)^\ha = M$.

\paragraph{} 
\label{tori-isog}
Write $T$ for the torus $\GB(A)^\ab = \GB(V')^\ab$.
Choose a faithful representation~$N$ of~$T$;
by the Tannakian formalism and \cref{hodge-is-motivated}
we may view~$N$ as an abelian \cm~motive over~$\CC$.

Since $\GB(W_F)$ is a torus, and $\tilde G$ is almost simple,
the representation
\[
 \GB(W_F) \times \tilde G \to \GL(V') = \GL(W_F \otimes_{F_S} V)
\]
has a finite kernel.
In other words, the natural map
$\GB(W_F) \times \tilde G \to \GB(A)$ is an isogeny,
and so is the natural map $\GB(W_F) \to T$.
We will need this fact in~\cref{same-ab-quot}.

\section{The main proposition} 

\begin{readme} 
 In this section we prove the main technical result of this paper.
 Its proof uses results of the preceding three sections.
\end{readme}

\begin{proposition} 
 \label{mtc-product-abelian-motives}
 Let $M_1$ and~$M_2$ be two
 geometrically irreducible hyperadjoint abelian motives
 over a finitely generated field $K \subset \CC$.
 Assume that $\MTC(M_1)$ and $\MTC(M_2)$ are true, and
 assume that $\Gl(M_1 \oplus M_2)$ is connected for all prime numbers~$\ell$.
 If there exists a prime number~$\ell$
 such that $\Gl(M_1 \oplus M_2) \subsetneq \Gl(M_1) \times \Gl(M_2)$,
 then $\HB(M_1) \cong \HB(M_2)$ as Hodge structures.
\end{proposition}

\paragraph{} 
\label{prp-strategy}
The proof of this \namecref{mtc-product-abelian-motives}
will take the remainder of this section.
Roughly speaking, the strategy is as follows:
\begin{enumerate}
 \item First we prove that $\End(M_1) = \End(M_2)$.
 \item Next we show that $\Hl(M_1) \cong \Hl(M_2)$ for all primes~$\ell$.
 \item We use this (and \cref{deldyn-local-global})
  to show that $M_1$ and~$M_2$
  have isomorphic Deligne--Dynkin diagrams over~$\QQ$.
 \item After that we run Deligne's construction (\cref{delignes-construction})
  on the motives~$M_i$;
  this leaves us with two complex abelian varieties~$A_1$ and~$A_2$
  such that $M_{i,\CC} = \HH^1(A_i)^\ha$.
 \item We replace $K$ be a finitely generated extension,
  such that $A_1$ and~$A_2$ are defined over~$K$.
 \item By carefully tracing the $\ell$-adic counterpart of the construction
  we show that $A_1$ and~$A_2$ have isomorphic $\ell$-adic Tate modules.
 \item Finally, we apply Faltings's results
  to deduce that $A_1$ and~$A_2$ are isogenous abelian varieties,
  which implies $\HB(M_1) \cong \HB(M_2)$.
\end{enumerate}
Sadly however, this strategy is slightly too optimistic.
It is not possible to work with the entire $\ell$-adic Galois representations:
we will have to focus our attention on a suitable summand.
This makes the proof quite technical.

\paragraph{} 
\label{isom-qcsr}
We first make some observations about~$M_1$ and~$M_2$.
For $i \in \{1,2\}$,
write $E_i$ for $\End(M_i)$, and
write $\Lambda_i$ for the set of finite places of~$E_i$.
Observe that $\HH_{\Lambda_i}(M_i)$ is a
quasi-compatible system of representations,
by \cref{abelian-motive-quasi-compatible-realisations}.

For a prime number~$\ell$, let $\Lambda_{i,\ell}$
denote the set of places~$\lambda \in \Lambda_i$ that lie above~$\ell$.
Recall that $\GB(M) = \Res_{E/\QQ}(G_i)$,
for some absolutely simple adjoint group~$G_i$ over~$E$.
Since $\MTC(M_i)$ holds
and Weil restriction of scalars commutes with base change,
we find that
\[
\Glc(M_i) = \GB(M) \otimes_\QQ \QQl
= \prod_{\lambda \in \Lambda_{i,\ell}}
\Res_{E_{i,\lambda}/\QQl} (G_i \otimes_E E_\lambda).
\]
In particular we have
$\Gl(\Hlambda(M_i)) = \Res_{E_{i,\lambda}/\QQl} (G_i \otimes_E E_\lambda)$
for every $\lambda \in \Lambda_{i,\ell}$.
If we write $\Glambda(M_i)$ for $\Glambda(\Hlambda(M_i))$
then the above computation implies that
$\Glambda(M_i)$ is equal to $G_i \otimes_E E_\lambda$,
an absolutely simple adjoint group over~$E_{i,\lambda}$.

Let $\ell$ be a prime number
such that $\Gl(M_1 \oplus M_2) \subsetneq \Gl(M_1) \times \Gl(M_2)$.
By Goursat's lemma (see \cref{Goursat}) this implies that there are
places $\lambda_1 \in \Lambda_{1,\ell}$ and $\lambda_2 \in \Lambda_{2,\ell}$
such that the projection of $\Gl(M_1 \oplus M_2)$ in
$\Res_{E_{1,\lambda_1}/\QQl} \GG_{\lambda_1}(M_1) \times
\Res_{E_{2,\lambda_2}/\QQl} \GG_{\lambda_2}(M_2)$
is the graph of an isomorphism of algebraic groups over~$\QQl$.
The derivative of this isomorphism at the identity element
is a $\QQl$-linear isomorphism
$\psi \colon \HH_{\lambda_1}(M_1) \to \HH_{\lambda_2}(M_2)$
of Galois representations of~$K$.
(Recall that $\HH_{\lambda_i}(M_i)$ is the adjoint representation
$\Lie(\Res_{E_{i,\lambda_i}/\QQl} \GG_{\lambda_i}(M_i))$.)

Note that $\psi$ induces an isomorphism
$f \colon \End_{\QQl}(\HH_{\lambda_1}(M_1))^{\Gamma_K} \to
\End_{\QQl}(\HH_{\lambda_2}(M_2))^{\Gamma_K}$.
Since these endomorphism algebras are commutative,
the isomorphism~$f$ does not depend on~$\psi$.
By \cref{recover-endomorphisms} we recover
$\lambda_i \colon E_i \into E_{i,\lambda_i} =
\End_{\QQl}(\HH_{\lambda_i}(M_i))^{\Gamma_K}$
as the subfield generated by coefficients of
characteristic polynomials of Frobenius elements
acting $E_{i,\lambda_i}$-linearly on~$\HH_{\lambda_i}(M_i)$.
Because $\psi$ commutes with the action of these Frobenius elements,
the map~$f$ restricts to a canonical isomorphism
$f \colon E_1 \to E_2$ that identifies $\lambda_1$ with~$\lambda_2$.
Write $E$ for $E_1 = E_2$,
and write $\lambda$ for $\lambda_1 = \lambda_2$.
We conclude that $\Hlambda(M_1) \cong \Hlambda(M_2)$
as $\lambda$-adic Galois representations.

Write $\Lambda$ for the set of finite places of~$E$.
We assumed that $\Gl(M_{1} \oplus M_{2})$
is connected for all prime numbers~$\ell$.
Because $\Hlambda(M_1)$ and $\Hlambda(M_2)$
are semisimple and quasi-compatible,
\cref{quasi-compatible-semisimple-isomorphic}
shows that $\HLambda(M_1)$ and $\HLambda(M_2)$
are isomorphic quasi-compatible systems of representations.
We have now completed the first two steps of the strategy
outlined in \cref{prp-strategy}.

\paragraph{} 
\label{isom-DD}
For $i = 1,2$, let $(\Delta_i,\mu_i)$
be the Deligne--Dynkin diagram associated with~$M_i$,
as in \cref{abmotdeldyn}.
Let $f \colon \pi_0(\Delta_1) \to \pi_0(\Delta_2)$ be the map
that is defined by the canonical identifications
$\pi_0(\Delta_i) \cong \Hom(E,\QQbar)$ of $\Gamma_\QQ$-sets.
Since $\HLambda(M_1)$ and $\HLambda(M_2)$ are isomorphic
as $E$-linear quasi-compatible systems of representations,
we get local isomorphisms
$\psi_\ell \in \Isom_f(\Delta_{1,\QQl},\Delta_{2,\QQl})^{\Gamma_{\QQl}}$.
By \cref{hodge-tate-cochar-deldyn} we have
$\psi_\ell(\mu_1) = \mu_2$ for all~$\ell$.
Thus we may apply \cref{deldyn-local-global}
to obtain an isomorphism
$\phi \in \Isom_f\big((\Delta_1,\mu_1),(\Delta_2,\mu_2)\big)^{\Gamma_\QQ}$,
such that $\phi_\lambda = \psi_\lambda$ for some finite place~$\lambda$ of~$E$.

\paragraph{} 
\label{run-del-constr}
Let $S_i$ be the subset of $\mu_i$-symplectic nodes of~$(\Delta_i,\mu_i)$.
We will now apply Deligne's construction (see \cref{delignes-construction})
to the motives~$M_1$ and~$M_2$.
To do so, we have to make three choices, as in \cref{choices}:
\begin{enumerate*}[label=(\textit{\roman*})]
\item choose a totally imaginary quadratic extension $F/E$, and
\item endow it with a partial \cm~type
 relative to $(\Delta_1,\mu_1) \stackrel{\phi}{=} (\Delta_2,\mu_2)$;
 finally,
\item choose the representations~$V_1$ and~$V_2$
 in such a way that they have the same dimension.
\end{enumerate*}
Let $F_{S_i}$ be an \'etale $E$-algebra such that
$\Hom(F_{S_i},\QQbar) \cong S_i$,
and write $W_{F,i}$ for $F \otimes_E F_{S_i}$.
Recall from \cref{WF} that $W_{F,i}$ carries a
fractional Hodge structure of weight~$1$ that is of \cm~type.
The construction produces two complex abelian varieties~$A_1$ and~$A_2$,
such that $\HB(A_i) = W_{F,i} \otimes_{F_{S_i}} V_i$,
and $M_{i,\CC} = \HH^1(A_i)^\ha$.

Replace~$K$ with a finitely generated extension of~$K$
such that $A_1$ and~$A_2$ are defined over~$K$.
We are now ready for the two final steps in~\cref{prp-strategy}.

\paragraph{} 
Let $\lambda$ be a finite place of~$E$ such that $\psi_\lambda = \phi_\lambda$,
and let $\ell$ be the residue characteristic of~$\lambda$.
Observe that $E$ acts on $M_i$,~$V_i$, $W_{F,i}$, and~$A_i$.
We have the following diagram:
\[
 \begin{tikzcd}
  && (W_{F_i} \otimes_E E_\lambda)^* \ar[d] \\
  && \Glambda(A_i)(E_\lambda) \ar[r,hook] \ar[dd,two heads]
  & \GL(\Hlambda(A_i)) \\
  \Gamma_K \ar[urr,dashed,bend left,"\tilde\rho_{i,\lambda}"]
  \ar[drr,bend right,"\rho_{i,\lambda}"] &
  \widetilde{\Glambda(M_i)}(E_\lambda) \ar[dr,two heads] \ar[ur] \\
  && \Glambda(M_i)(E_\lambda) \ar[r,hook] & \GL(\Hlambda(M_i))
 \end{tikzcd}
\]
If we temporarily forget the dashed arrow $\tilde\rho_{i,\lambda}$,
then we can identify the two diagrams for $i = 1$ and $i = 2$,
to obtain a diagram:
\[
 \begin{tikzcd}
  && W_{F,\lambda}^* \ar[d] \\
  && \Glambda'(E_\lambda) \ar[r,hook] \ar[dd,two heads] & \GL(\Hlambda') \\
  \Gamma_K \ar[drr,bend right,"\rho_\lambda"] &
  \widetilde{\Glambda}(E_\lambda) \ar[dr,two heads] \ar[ur] \\
  && \Glambda(E_\lambda) \ar[r,hook] & \GL(\Hlambda)
 \end{tikzcd}
\]
These identifications are justified as follows:
\begin{itemize}
 \item In \cref{isom-qcsr} we showed that there is an
  isomorphism $\Hlambda(M_1) \to \Hlambda(M_2)$
  of $\lambda$-adic Galois representations,
  identifying $\Glambda(M_1)$ with~$\Glambda(M_2)$
  and $\rho_{1,\lambda}$ with~$\rho_{2,\lambda}$.
  We may also canonically identify
  $\widetilde{\Glambda(M_1)}$ and $\widetilde{\Glambda(M_2)}$.
  Call the resulting objects
  $\Hlambda$, $\Glambda$, $\rho_\lambda$, and~$\widetilde{\Glambda}$
  respectively.
 \item This isomorphism also gives an isomorphism
  of Deligne--Dynkin diagrams
  $\psi_\lambda \colon (\Delta_1,\mu_1)_\lambda \to (\Delta_2, \mu_2)_\lambda$
  over~$\QQl$. In \cref{isom-DD} we showed that
  $\psi_\lambda$ extends to an isomorphism
  of Deligne--Dynkin diagrams
  $\phi \colon (\Delta_1,\mu_1) \to (\Delta_2, \mu_2)$ over~$\QQ$.
 \item This allows us to identify
  the $E$-algebras $W_{F,1}$ and~$W_{F,2}$.
  Call the resulting algebra~$W_F$
  and write $W_{F,\lambda}$ for $W_F \otimes_E E_\lambda$.
 \item By Deligne's construction
  we can identify $\Hlambda(A_1)$ with $\Hlambda(A_2)$
  as representations of~$\widetilde{\Glambda}$
  with an action of~$W_{F,\lambda}$.
  Call this representation~$\Hlambda'$.
 \item Finally, $\Glambda(A_i)$ is the Zariski closure
  of the image of~$W_{F,\lambda}^*$ and ~$\widetilde{\Glambda}$
  in $\GL(\Hlambda(A_i))$.
  Thus we may identify $\Glambda(A_1)$ with $\Glambda(A_2)$
  and call the resulting group $\Glambda'$.
\end{itemize}
The dashed arrows $\tilde\rho_{1,\lambda}$
and~$\tilde\rho_{2,\lambda}$ in the original diagrams
can be identified with maps $\Gamma_K \to \Glambda'(\QQl)$
that lift~$\rho_\lambda$.
It is our goal to show that they are identical.

\paragraph{} 
To show that $\tilde\rho_{1,\lambda}$ and~$\tilde\rho_{2,\lambda}$
are identical,
we employ the following strategy:
Recall that $\Glambda'$ is a finite cover of
$(\Glambda')^\ad \times (\Glambda')^\ab$.
Since $\tilde\rho_{1,\lambda}$ and~$\tilde\rho_{2,\lambda}$
lift the map~$\rho_\lambda$
we know that the their composition with the canonical quotient map
$\Glambda' \to (\Glambda')^\ad$ is independent of $i = 1,2$.
If we also show that their composition with the canonical quotient map
$\Glambda' \to (\Glambda')^\ab$ is independent of $i = 1,2$,
then we can show that
$\tilde\rho_{1,\lambda}$ and~$\tilde\rho_{2,\lambda}$
are identical,
at least after replacing $K$ by a finite field extension of~$K$.
In actually executing this strategy, we play a minor variation.

\paragraph{} 
\label{same-ab-quot}
Recall from \cref{tori-isog}
that for $i = 1,2$,
the natural map $p_i \colon \GB(W_F) \to \GB(A_i)^\ab$ is an isogeny.
Therefore there exists a torus~$T$,
together with isogenies $q_i \colon \GB(A_i)^\ab \to T$
such that $q_1 \circ p_1 = q_2 \circ p_2$.
Let $N$ be a faithful representation of~$T$.
Since $N$ is a representation of~$\GB(W_F)$
we may canonically view $N$ as a fractional pre-Hodge structure,
via the Tannakian formalism.
However, since $N$ is a representation of $\GB(A_i)$
we know that $N$ is a (classical) Hodge structure,
and this Hodge structure is independent of~$i$,
because $q_1 \circ p_1 = q_2 \circ p_2$. 
Note that $N$ is a \cm~Hodge structure, by definition.
Finally, by \cref{hodge-is-motivated} we view $N$ as a complex \cm~motive.

Replace~$K$ with a finitely generated extension of~$K$
such that the motive~$N$ is defined over~$K$.
Observe that the composite map
\[
 \Gamma_K \stackrel{\tilde\rho_{i,\ell}}{\longrightarrow}
 \underbrace{\Gl(A_i)(\QQl)}_{\GB(A_i)(\QQl)} \longrightarrow
 \underbrace{\Gl(A_i)^\ab(\QQl)}_{\GB(A_i)^\ab(\QQl)}
 \stackrel{q_{i,\ell}}{\longrightarrow}
 \underbrace{\Gl(N)(\QQl)}_{T(\QQl)}
\]
defines the Galois representation on~$\Hl(N)$
and therefore is independent of~$i = 1,2$.

Recall that we have quotient maps
$\Gl(A_i) \onto \Res_{E_\lambda/\QQl}\Glambda(A_i)
= \Res_{E_\lambda/\QQl}\Glambda'$,
and therefore
$\Gl(A_i)^\ab \onto (\Res_{E_\lambda/\QQl}\Glambda')^\ab$.
We claim that there exists a finite quotient
$(\Res_{E_\lambda/\QQl}\Glambda')^\ab \onto T'$
such that the composite map
\[
 \Gamma_K \stackrel{\tilde\rho_{i,\ell}}{\longrightarrow}
 (\Res_{E/\lambda/\QQl}\Glambda')(\QQl) \longrightarrow
 (\Res_{E_\lambda/\QQl}\Glambda')^\ab(\QQl)
 \longrightarrow T'(\QQl)
\]
is independent of $i = 1,2$.
Indeed, one may construct $T'$ as follows:
Let $K_i$ be the kernel of the isogeny
$q_{i,\ell} \colon \Gl(A_i)^\ab \onto \Gl(N)$,
and let $K_i'$ denote its image under the quotient map
$\Gl(A_i)^\ab \onto (\Res_{E_\lambda/\QQl}\Glambda')^\ab$.
Then $K_1' \cdot K_2'$ is a finite subgroup of
$(\Res_{E_\lambda/\QQl}\Glambda')^\ab$,
and we define $T'$ as the quotient.
By construction there is a map $\Gl(N)(\QQl) \to T'(\QQl)$
and the claim follows.

\paragraph{} 
We are almost done!
Write $\pi^\ab$ for the map
$\Glambda'(E_\lambda) = (\Res_{E_\lambda/\QQl}\Glambda')(\QQl) \to T'(\QQl)$
and $\pi^\ad$ for the map
$\Glambda'(E_\lambda) \to (\Glambda')^\ad(E_\lambda) = \Glambda(E_\lambda)$.
We have now proven that the compositions
$\pi^\ab \circ \tilde\rho_{i,\lambda}$
and $\pi^\ad \circ \tilde\rho_{i,\lambda}$
do not depend on $i = 1,2$.

Recall that $T'$ is a quotient of~$(\Res_{E_\lambda/\QQl}G')^\ab$
by a finite group.
Hence the group $\alpha = \ker(\pi^\ad) \cap \ker(\pi^\ab)$ is
a finite abelian subgroup of the centre of~$\Glambda'(E_\lambda)$.
Define $\xi \colon \Gamma_K \to \Glambda'(E_\lambda)$
via $\xi(g) = \tilde\rho_{1,\lambda}(g) \cdot \tilde\rho_{2,\lambda}(g)^{-1}$.
This map $\xi$ takes values in~$\alpha$,
and is a homomorphism since $\alpha$
is contained in the centre of~$\Glambda'(E_\lambda)$.
We conclude that after replacing $K$ by a finite extension,
the homomorphisms $\rho_{1,\lambda}$ and~$\rho_{2,\lambda}$ are the same.

This means that $\Hlambda(A_1)$ and~$\Hlambda(A_2)$
are isomorphic as $\lambda$-adic Galois representations.
Another application of \cref{quasi-compatible-semisimple-isomorphic}
shows that $\HLambda(A_1)$ and~$\HLambda(A_2)$
are isomorphic $E$-rational quasi-compatible systems of representations.
In particulare $\Hl(A_1)$ and~$\Hl(A_2)$ are isomorphic Galois representations.
Finally, Faltings's theorem
(Korollar~2 in~\S5 of~\cite{Fal83}, see also~\cite{Fal84})
implies that $A_1$ and~$A_2$ are isogenous abelian varieties.
Since $M_i = \HH(A_i)^\ha$,
we conclude that $\HB(M_1)$ and~$\HB(M_2)$ are isomorphic Hodge structures.
This completes the proof of \cref{mtc-product-abelian-motives}.

\section{The main theorem} 

\begin{lemma} 
 \label{subgroup-surjective-projection-binary-products}
 Let $K$ be a field of characteristic~$0$.
 For $i = 1, \dots, n$, let $G_i$ be a simple linear algebraic group over~$K$.
 Let $G$ be a subgroup of $G_1 \times \cdots \times G_n$,
 such that $G$ surjects onto $G_i$ for $1 \le i \le n$
 and $G$ surjects onto $G_i \times G_j$ for $1 \le i < j \le n$.
 Then $G = G_1 \times \cdots \times G_n$.
 \begin{proof}
  It suffices to prove the analogous statement for Lie algebras.
  This is precisely
  the lemma in step~3 on pages~790--791 of~\cite{Ribet}.
 \end{proof}
\end{lemma}

\begin{lemma} 
	\label{mtc-finite-product-abelian-motives}
	Let $K \subset \CC$ be a finitely generated field.
	Let $M_{i}$, with $i \in I$, be a finite collection of
	irreducible hyperadjoint abelian motives over~$K$.
	Write $M = \bigoplus M_{i}$.
	If $\MTC(M_{i})$ is true for all $i \in I$,
	then $\MTC(M)$ is true.
	\begin{proof}
		By replacing $K$ with a finite extension,
		we may and do assume that $M_i$ is \emph{geometrically} irreducible
		for all $i \in I$.
		We also may and do assume that for $i,j \in I$
		we have $M_{i} \cong M_{j} \iff i = j$.
		Recall that there is a natural injection
		$\Gl(M) \into \prod_{i \in I} \Gl(M_{i})$
		and the image projects surjectively onto the factors $\Gl(M_{i})$.
		By \cref{mtc-product-abelian-motives},
		we know that if $i,j \in I$ are two different indices,
		then $\Gl(M_{i} \oplus M_{j}) \cong \Gl(M_{i}) \times \Gl(M_{j})$;
		in other words,
		$\Gl(M)$ surjects onto $\Gl(M_{i}) \times \Gl(M_{j})$.
		By \cref{subgroup-surjective-projection-binary-products}
		we conclude that $\Gl(M) \cong \prod_{i \in I} \Gl(M_{i})$.
	\end{proof}
\end{lemma}

\begin{theorem} 
	\label{mtc-abelian-motives-tannakian-subcategory}
	Fix a finitely generated field $K \subset \CC$.
	Let $\mathcal{M}_K \subset \Mot_K$ be the full subcategory
	of abelian motives for which the Mumford--Tate conjecture is true.
	Then the category~$\mathcal{M}_K$ is a Tannakian subcategory of $\Mot_K$.
	\begin{proof}
		The category~$\Mot_K$ is semisimple,
		so subquotients are direct summands.
		It is clear that the subcategory~$\mathcal{M}_K$
		is closed under duals, tensor powers, and direct summands.
		Let $M_1$ and~$M_2$ be two objects in~$\mathcal{M}_K$.
		We need to show that $M_1 \oplus M_2$
		and $M_1 \otimes M_2$ are objects in $\mathcal{M}_K$.
		Observe that $M_1 \otimes M_2$ is a direct summand of
		$(M_1 \oplus M_2)^{\otimes 2}$.
		Thus we are done if we show that the Mumford--Tate conjecture is true for
		$M = M_1 \oplus M_2$.
		By \cref{mtc-adjoint-motive} it suffices to prove $\MTC(M^\ha)$.
		Decompose $M^\ha = \bigoplus_{i \in I} M'_i$
  into a sum of irreducible motives.
  Note that for every $i \in I$, the motive~$M'_i$ is hyperadjoint,
  and $\MTC(M'_i)$ holds
  since the irreducible motive~$M'_i$ is a summand of $M_1^\ha$ or~$M_2^\ha$
  by \cref{ha-summands}.
		Now the result follows from \cref{mtc-finite-product-abelian-motives}.
	\end{proof}
\end{theorem}

\begin{remark} 
 We give a short and incomplete list of examples of varieties
 for which it is known that their motives
 are objects of the category~$\mathcal{M}_K$
 of the preceding \namecref{mtc-abelian-motives-tannakian-subcategory}.
 (See \cite{Mo17} for a recent overview of the state of the art.)
 \begin{enumerate}
  \item Simple abelian varieties of prime dimension (\cite{Ta83}).
  \item Abelian varieties of dimension~$g$ with trivial endomorphism ring
   such that $2g$ is neither a $k$-th power for some odd $k > 1$,
   nor of the form $\binom{2k}{k}$ for some odd $k > 1$
   (theorem~5.14 of~\cite{Pi98}).
  \item K3~surfaces (\cite{MTCK3I} and~\cite{MTCK3II}).
  \item Cubic fourfolds (theorem~1.6.1 of~\cite{An96}).
 \end{enumerate}
 We remind the reader that the Mumford--Tate conjecture
 for an abelian fourfold~$A$ with trivial endomorphism ring
 is still an open problem.
\end{remark}

\begin{theorem} 
 \label{mtcaxa}
 Let $K$ be a finitely generated subfield of~$\CC$.
 Let $A_i$, $i \in I$, be a finite collection of abelian varieties over~$K$,
	and write $A = \prod_{i \in I} A_i$.
 Assume that $\MTC(A_i)$ is true for all $i \in I$.
 Then $\MTC(A)$ is also true.
 \begin{proof}
  This is an immediate consequence of
  \cref{mtc-abelian-motives-tannakian-subcategory}.
 \end{proof}
\end{theorem}

\printbibliography

\end{document}